\documentclass[square,numbers]{oupau}
\usepackage{amsmath, amsfonts,amssymb, latexsym}
\usepackage[all]{xy}
\usepackage{color} 
\usepackage{graphicx} 
\usepackage{graphics}

\newtheorem{defn}[theorem]{Definition}
\newtheorem{thm}[theorem]{Theorem}

\newcommand{\ag}{{\alpha}}
\newcommand{\al}{{\alpha}}

\newcommand{\be}{{\beta}}
\newcommand{\bg}{{\beta}}

\newcommand{\om}{{\omega}}
\newcommand{\og}{{\omega}}
\newcommand{\io}{{\iota}}

\newcommand{\eps}{{\varepsilon}}

\newcommand{\dg}{{\delta}}
\newcommand{\eg}{{\epsilon}}
\newcommand{\ga}{{\gamma}}

\newcommand{\la}{{\lambda}}
\newcommand{\La}{{\Lambda}}

\newcommand{\sg}{{\sigma}}

\newcommand{\iso}{{\cong}}

\newcommand{\Cc}{{\mathcal C}}
\newcommand{\Dd}{{\mathcal D}}
\newcommand{\Ll}{{\mathcal L}}

\newcommand{\Mm}{{\mathcal M}}

\newcommand{\Ss}{{\mathcal S}}

\def	\half	{\frac{1}{2}}

\newcommand{\bb}[1]{\ensuremath{\mathbb{#1}}}

\renewcommand{\to}{\longrightarrow}

\newcommand{\PD}{{\rm PD} }

\newcommand{\CP}[1]{\bb{C}P^{#1}}
\newcommand{\Symp}{\mathrm{Symp}}
\newcommand{\Diff}{\mathrm{Diff}}
\newcommand{\red}[2]{\overline{#1}_{#2} }
\newcommand{\Bl}[2]{{Bl}_{#2}(#1)}

\newcommand{\hs}[1][2n]{\ensuremath{\mathbb{H}\mathrm{Symp}_{#1}}}
\newcommand{\HS}[1][2n]{\ensuremath{\mathbf{H}\mathrm{Symp}_{#1}}}
\newcommand{\qu}{/\kern-.7ex/}

\begin{document}
\title{Classifying semi-free Hamiltonian $S^1$-manifolds}
\author{Eduardo Gonz\'alez}
\address{Mathematics Department, UMASS Boston.
}
\correspdetails{eduardo@math.umb.edu}
\received{1 January 2009}
\revised{11 January 2009}
\accepted{21 January 2009}


\date{January, 2009}
\begin{abstract}
  In this paper we describe a method to establish when a symplectic
  manifold $M$ with semi-free Hamiltonian $S^{1}$-action is unique
  up to isomorphism (equivariant symplectomorphism).  This will rely
  on a study of the symplectic topology of the reduced spaces. We
  prove that if the reduced spaces satisfy a rigidity condition, then
  the manifold $M$ is uniquely determined by fixed point data.  In
  particular we can prove that there is a unique family up to
  isomorphism of 6-dimensional symplectic manifolds with semi-free
  $S^1$-action and isolated fixed points.
\end{abstract}

\maketitle

\section{Introduction.}\label{s:intro}

Let $(M,\om)$ be a $2n$ dimensional compact, connected, symplectic
manifold, and let $\{\La_t\}=\La: S^1\to \mathrm{Symp}(M,\om)$ denote
a symplectic circle action on $M$, that is, if $X$ is the vector field
generating the action, then $\mathcal{L}_X \om = d \io_X \om = 0$.
Recall that the action is \textbf{semi-free} if it is free on $M
\backslash M^{S^1}$, i.e.  the only {\em weights} at every fixed point
are $\pm 1, 0$. A circle action is said to be Hamiltonian if there is
a $C^{\infty}$ function $H:M \to \bb{R} $ such that $\io_X \om = -d
H$.  Such a function is called a Hamiltonian for the action.  $H$ is
unique if we normalize it by requiring its minimum to be $0$.  $H$ is
a perfect Morse-Bott function and its critical points are the fixed
points of $\Lambda$. Let $\Cc(M)$ denote the critical values of
$H$. This set is determined by the cohomology class of $\omega$
(cf. Lemma \ref{l:criticalv}).  We will denote by $\HS$ the class of
manifolds $(M,H,\om)$ with semi-free Hamiltonian circle actions with
normalized Hamiltonians.

In this paper we will provide a mechanism that classifies, under
restricted conditions, Hamiltonian $S^1$-manifolds $(M,H,\om)$ up to
\textbf{isomorphism}, i.e. equivariant symplectomorphism. The idea is
to investigate how to reconstruct $M$ from its \textsl{local data}
$\Ll(M)$ and what type of information determines the local
data. $\Ll(M)$ can be thought as an atlas for Hamiltonian
$S^1$-actions. It will be, roughly speaking, given by neighborhoods of
the critical levels of $H$ called \textsl{critical germs} and open
submanifolds of regular points called \textsl{slices}. Precise
definitions can be found in section \S\ref{s:general}. Our first
result can be stated roughly as follows (see Theorem \ref{t:glu}).

\begin{theorem}
  The local data $\Ll(M)$ determines the Hamiltonian $S^1$-manifold
  $(M,H,\omega)$ up to isomorphism.
\end{theorem}

Although we carry out the proof of this theorem in the semi-free case,
it holds without this assumption.

Our second main objective is to describe the local data in terms of
more intrinsic information from the critical levels. To accomplish
that, we have to make certain concessions. First, to describe the
germs we restrict to the case when all non-extremal critical levels
contain fixed-point components of Morse (co)index at most 2.  In this
case, as shown in \S\ref{sss:smooth}, all the reduced spaces
$\red{M}{\la}$ at critical values $\lambda$ have a natural structure
of a smooth symplectic manifold $(\red{M}{\lambda},
\overline{\omega}_\lambda)$.  For clarity we will further assume that
the critical levels are \textbf{simple} (see Section
\ref{ss:non-simple} for the general case), this means that all the
fixed point components in the same level have a common index
$i_\lambda$. Under this assumption, the diffeomorphism type of the
reduced spaces $\red{M}{t}$ is constant on the semi-closed interval
$(\lambda-\eps,\lambda]$. In general, even in the non-simple case, the
diffeomorphism type of the reduction bundle $P_t:=H^{-1}(t)\to
H^{-1}(t)/S^1=\red{M}{t}$ is constant for $t$ in any interval of
regular values of $H$. We denote by $P_{\lambda}^{-}$ and $P_{\lambda}^{+}$
the diffeomorphism types for regular values immediately below and
above $\lambda$, respectively. We let $e(P_\lambda^{\pm})\in
H^2(\red{M}{\lambda})$ denote the Euler class of $P_\lambda^{\pm}$.

Before we describe the isomorphism type of the germs we need the
following definition. 







\begin{defn}\label{d:fpd-s}
  Let $(M,H,\om)$ denote a Hamiltonian $S^1$-manifold. Suppose that
  all the non-extremal critical levels are simple and of index 2. If
  $\lambda$ denote a non-extremal critical value of $H$, the
  \textbf{fixed point data of $M$ at $\lambda$} consists of the
  following.

  \begin{itemize}
  \item The critical value $\lambda$ and its index $i_\lambda$.
  \item For each fixed component $F$ at $\lambda$, the tuple
    $(\red{M}{\la}, \overline{F}, \overline{\om}_\lambda )$ of the
    reduced space and the image of the fixed point component embedded
    as a symplectic submanifold. (This information contains the normal
    bundle $(N_{ \overline{F} / {\red{M}{\la}}})$.)
  \item The bundle $P_\lambda^{-}$.
\end{itemize}

A similar definition applies if $\lambda$ has coindex 2.  For extremal
fixed point components, let $\rm{ext}$ denote either $\rm{max}$ or
$\rm{min}$. One has that $\red{M}{\la_{\rm{ext}}} =
\red{F}{\la_{\rm{ext}}}$. Then, we take the fixed point data as
\begin{itemize}
\item The critical value $\lambda_{\rm{ext}}$.
\item The symplectic manifold
  $(\red{F}{\la_{\rm{ext}}},\overline{\om}_{\la_{\rm{ext}}})$.
\item The symplectic normal bundle $N_{F_{\la_{\rm{ext}}}/ M}$ of
  $F_{\la_{\rm{ext}}}$ in $M$.
\end{itemize}
In this situation the index $i_{\lambda_{\rm{ext}}}$ can be recovered
by the normal bundle $N_{F_{\la_{\rm{ext}}}/ M}$, and thus it is
not necessary.  The \textbf{fixed point data} of $M$ consists of the
collection of all of the above tuples for $\la\in \Cc(M)$.
\end{defn}

It is important to recall that in general the fixed point set at
non-extremal critical levels might be disconnected and in the
non-simple case they might have different Morse indexes. The extremal
fixed point sets are always connected.

It is already stated in Guillemin-Sternberg's \cite{GS1} that the
isomorphism type of the critical germs is determined by the fixed
point data. As we will see, under further assumptions, we can classify
manifolds knowing less information in the fixed point data.  As
follows.

\begin{defn}\label{d:mfpd}
  The \textbf{small fixed point data} of $M$ consist of
  the following information at each  $\lambda$.
  \begin{itemize}
  \item If $\lambda$ is the minimum, the same information as above.
  \item If $\lambda$ is non-extremal, then for each fixed point component
    $\overline{F}$ at $\lambda$ we consider the following data.
    \begin{itemize}
    \item If $\dim \overline{F}>0$, then we keep the same information
      above but excluding $e(P_\lambda)$.
    \item If $\dim \overline{F}=0$ then we only consider $\lambda$ and
      its index $i_\lambda$. 
    \end{itemize}
  \item If $\lambda$ is the maximum, then we only take $\lambda$ and
    the symplectic manifold $(\overline{F}_{\rm
      max},\overline{\omega}_{\rm max})$. 
  \end{itemize}
\end{defn}

To understand the regular slices, we first recall the content of
Duistermaat-Heckman theorem. On an interval $I=[t_0,t_1]$ of regular
values of $H$, all the reduced bundles $P_t\to \red{M}{t}, t\in I$ are
diffeomorphic, say to a fixed one $P\to B$. Moreover the symplectic
structure $\omega$ on $M$ yields a family $\{\overline{\omega}_{t\in
  I}\}$ of reduced symplectic forms on $B$. The pair
$(B,\{\overline{\omega}_t\})$ and the bundle $P\to B$ determines the
isomorphism class of the regular slice $H^{-1}(I)$. Theorem
\ref{p:equiv} extends this result, showing which families of
symplectic forms give isomorphic regular slices in the case when one
has the following property.


\begin{definition}\label{d:rigid}
  Let $B$ be a manifold and $\{\overline{\omega}_t\}$ be a family of
  symplectic structures on $B$ smoothly parametrised by $t\in I$. The
  pair $(B,\{\overline{\omega}_t\})$ is said to be \textbf{rigid} if
  \begin{enumerate}
  \item\label{d:pathconnected} $\Symp(B,\overline{\omega}_t)\cap {\rm
      Diff}_0(B)$ is path connected for all $t\in I$.
  \item\label{d:deformation} Any deformation between any two
    cohomologous symplectic forms which are deformation equivalent to
    $\omega_{t_0}$ on $B$ may be homotoped through deformations with
    fixed endpoints into an isotopy.
  \end{enumerate}
\end{definition}

Rigidity can be understood as the symplectic analogue of complex
rigidity, in the sense that even if we deform the symplectic forms, we
would obtain isotopic forms. Property \eqref{d:pathconnected} is in
fact non-generic, and Seidel \cite{Se-Le} has shown it is not
satisfied for the monotone blow up $\CP{2}\#5\overline{\CP{2}}$.

Once we have understood the regular slices and the germs, we then glue
them near the critical levels, a construction that is already used by
Hui Li \cite{Li2}. Assuming all the reduced spaces are rigid, the
resulting manifold is independent of the gluing information and we
prove the following theorem.

\begin{thm}[Weak classification Theorem]\label{t:intro1}
  Let $(M,H,\om)\in \HS$.  Suppose that all critical levels at
  $\lambda\in\Cc(M)$ are simple. Suppose further that for any two
  consecutive $\lambda',\lambda\in \Cc(M)$, $\lambda'<\lambda$, the
  pair $(\red{M}{\la},\{\overline{\omega}_t\}_{t\in (\lambda',\lambda)})$
  is rigid. Then $(M,H,\om)$ is determined by its fixed point data up
  to equivariant symplectomorphism.
\end{thm}



In the case when $\dim M=6$ all the non-extremal fixed-point
components have (co)index 2. Following the notation of Karshon-Tolman
\cite{KaTo-Ce}, this corresponds to complexity two (half-dimension of
the reduced spaces) Hamiltonian symplectic manifolds. Therefore we can
use established results on the symplectic structures
and on the topology of the group of symplectomorphisms. 
When the fixed point components are simple enough it is possible to
describe some of the fixed point data at each level in terms of fixed
point data associated to the previous critical level. Therefore the
isomorphism type only depends on the small fixed point data of
Definition \ref{d:mfpd}. We have got the following result.

\begin{theorem}\label{t:baby2}
  Let $(M,H,\omega)\in\HS[6]$. Suppose that all critical levels at
  $\lambda\in\Cc(M)$ are simple and that for any two consecutive
  $\lambda',\lambda\in \Cc(M)$, $\lambda'<\lambda$, the pair
  $(\red{M}{\la},\{\overline{\omega}_t\}_{t\in (\lambda',\lambda)})$
  is rigid. If the fixed point sets $(\red{F}{\la},
  \overline{\om}_{F_\la})$ are either surfaces or isolated fixed
  points, then the isomorphism class of $(M,H,\omega)$ is uniquely
  determined by the small fixed point data.
\end{theorem}


In \cite{Li} Hui Li analysed manifolds using even weaker fixed point
data, for instance, the normal bundle of all the fixed point sets is
not considered. Her information is in spirit purely topological and in
that case she has constructed manifolds with the same fixed point data
(in her sense) but not diffeomorphic, she calls this phenomena a
\emph{twist}. Theorem \ref{t:baby2} shows that provided all reduced
spaces are rigid and we know the fixed point data as in Definition
\ref{d:fpd-s}, such twist cannot exist.

We will push our method even further by classifying a family of
manifolds with minimal fixed point information. Let $Y^{n}=S^{2}\times
\dots \times S^{2}$ be the $n$-fold product of spheres and let $\sg$
be the canonical area form on $S^{2}$. Provide $Y^{n}$ with the
product symplectic form $\la_{1}\sg\times \dots \times \la_{n}\sg$
that takes the value $\la_{i}>0$ on each of the spheres of $Y^{n}$.
Let the circle act diagonally on $Y^{n}$ in the standard semi-free
Hamiltonian fashion. $Y^{n}$ is the only known example of a
$2n$-dimensional symplectic manifold that admits a semi-free circle
action with isolated fixed points. Thus, it is natural to ask if this
is the only manifold up to equivariant symplectomorphisms that has
this property. In the case $n=2$ the methods of Karshon \cite{K}
answer this question positively. In the present paper we establish the
result for $n=3$.

\begin{theorem}\label{t:intro}
  Let $M$ be a 6-dimensional symplectic manifold with a semi-free
  circle action that has isolated fixed points. Then, $M$ is
  equivariantly symplectomorphic to $Y^{3}$ with the canonical product
  form for some $\la_{i}$.
\end{theorem}

Plenty of information was known about these manifolds.  In \cite{TW}
it is proved that under the hypothesis of Theorem \ref{t:intro} the
action must be Hamiltonian. The $\la_i$ are in fact the critical
values of the Hamiltonian function of the (only) three fixed points of
index 2, assuming the minimum to be zero.  Moreover, there is an
isomorphism of equivariant cohomology rings $H_{S^{1}}^*(M) \cong
H_{S^{1}}^*(Y^3)$. The isomorphism is such that takes Chern classes
into Chern classes, thus by Wall \cite{W} $M$ is diffeomorphic to
$Y^3$. More recently and in higher dimensions, a result of Ilinski
reproved by Matsuda and Panov \cite[Corollary 4.9]{MaPa-Se06} shows
that if $M$ is known to be a $2n$ dimensional toric space with a
circle subgroup acting semi-freely and with isolated fixed points,
then $M$ must be diffeomorphic to $Y^n$. Nothing was known about the
symplectomorphism type of $M$.  An approach to understand its
symplectomorphism type was given in \cite{G1}, where we proved that
the quantum cohomology ring of $M$ and $Y^3$ are isomorphic.

The proof of Theorems \ref{t:intro1} and \ref{t:intro} will be done in
stages. We will show that one can construct the desired isomorphism
starting from the minimum and extending it by gluing regular slices
and germs. As we mention above, this local information only depends on
the fixed point data and the symplectomorphism type of the reduced
manifolds, which in this very particular example, are $\CP{2}$ and its
blow-ups, whose symplectic topology is well known. As we will see,
since the fixed components are isolated points, the local data at a
critical level will depend only on the data of the previous one. Thus
we can ``bootstrap'' the construction to extend the isomorphism by
attaching a symplectic manifold whose symplectomorphism type is
determined by previous information.  For instance, one knows that
after passing a critical level, the reduced manifolds are related by a
blow-up and blow-down, and thus the resulting manifold after passing a
critical level is determined by the previous data.

The results in \cite{TW},\cite{G1} work in all dimensions.  However,
we do not expect the techniques presented in the current paper to
extend easily to higher dimensions.  This is because our rigidity
arguments strongly rely on results in four dimensional symplectic
geometry, which as of now, do not have higher dimensional versions.
We finally note that the proof of Theorem \ref{t:intro} does not use
Wall's theorem, in contrast to previous methods.

As a last remark, we point out that it is possible to remove the
semi-free assumption from our results. In this case one has to deal
with generalizations of our tools to the orbifold category, as the
work of L. Godinho \cite{God1,God3} and W. Chen \cite{Ch-Ps}. This has
already been explored by the recent work of McDuff \cite{Mc-So09} in
some particular case.

\vspace{.1in}

\noindent {\bf Acknowledgments:} I thank Dusa McDuff and Chris
Woodward for encouragement and support during the preparation and
revision of this paper.  I also owe a special debt of gratitude to Sue
Tolman for her illuminating objections and in-depth revision of the
first version of this paper.  Finally, I thank Hui Li and Martin
Pinsonnault for discussions on an early draft and to the anonymous
Referee for pointing out many inconsistencies in the manuscript.


\section{General Setting.}\label{s:general}

In this section we will be using the term Hamiltonian $S^1$-manifold
for triples $(M,H,\om)$, where $M$ is a closed, smooth,
connected $2n$-manifold, $\om$ a symplectic form on $M$ and $H$ is a
normalized Hamiltonian function on $M$ that generates a circle action
compatible with $\om$ on $M$.  We will assume the action to be
semi-free, just to be consistent with the main objective of the paper.
For this section this hypothesis can be removed.  Although we will be
working with general, not necessary closed manifolds, we will always
think of them as (isomorphic to) submanifolds of a closed one, say
$M$.  Therefore, we sometimes will not make explicit the presence of
the symplectic form.  We will use the notation $\hs$ referring to the
class of Hamiltonian $S^{1}$-manifolds, closed or not, up to
isomorphism.  Here $(M,H,\om)$ and $(M',H',\om')$ are
\textbf{isomorphic} if there exist an equivariant diffeomorphism
$f:M\to M'$ such that $f^*(\om')=\om$.  Note that $f$ is equivariant
if and only if $H\circ f=H'$.

Let $(M,L,\om)$ be a closed Hamiltonian $S^1$-manifold and denote by
\[\Cc(M)=\{ 0=\la_0 < \dots < \la_s \}\] the collection of critical
values of $L$. Thus $L(M)=[\la_0,\la_s]$.

\begin{lemma}\label{l:criticalv}
  The set $\Cc(M)$ is invariant under isomorphism and is determined by
  the cohomology class of the symplectic structure.
\end{lemma}

\begin{proof}
  The first statement is obvious. For the second, choose any invariant
  metric in $M$. Let $p\in M^{S^1}$ a fixed point and take any
  gradient line of $L$ joining $p$ to a point $q$ in the minimum. Let
  $S$ denote the 2-cycle (sphere) obtained by rotating the gradient
  line by the circle action. By the Hamiltonian assumption its
  symplectic area is given by
  \[
  \int_{S}\omega=L(p)-L(q)=L(p).
  \]
  This determines $L(p)$ purely on cohomological data. For more
  details the reader can consult \cite{G1}.
\end{proof}

We now describe how to get pieces of $M$ by \textsl{localising to
  critical levels using the Hamiltonian}. This is, for each $\eps>0$
consider the neighborhood $L^{-1}(\la-\eps,\la+\eps)$ of the level set
$L^{-1}(\la)$, and for two consecutive values $\la,\la'\in \Cc(M)$
take the open submanifold $L^{-1}(\la,\la')$. In this paper we are
interested on knowing what is needed to reconstruct $M$ if one just
knows the isomorphism type of the pieces $L^{-1}(\la-\eps,\la+\eps)$
and $L^{-1}(\la,\la')$. In general one would need to specify how to
glue these pieces to get back to $M$. To see what type of gluing maps
are allowed, suppose that $(Y,H),(Z, K)$ are manifolds isomorphic to
the pieces $L^{-1}(\la-\eps,\la+\eps)$ and $L^{-1}(\la,\la')$
respectively. Note that one is tempted to glue $Y$ and $Z$ along the
overlap, that is we want to identify $H^{-1}(\la,\la+\eps)$
and $K^{-1}(\la,\la+\eps)$ through isomorphisms, to obtain back a
manifold isomorphic to $L^{-1}(\la-\eps,\la')$ such that $Y$ and $K$
are symplectic submanifolds. Therefore one needs to consider a maximal
set of gluing maps with this property. One may think of these data as
a \textsl{symplectic atlas of compatible Hamiltonian charts} for
$M$. We now provide the precise formalism that allows this to work.

\begin{defn}\label{d:1}
Let   $\la\in \bb{R}$ and let $\eps_0>0$. 
\begin{itemize}
\item[(i)] A \textbf{cobordism} at $\la$ is a tuple $(Y,H, \eg)$ such
  that $0<\eg<\eps_0$ and $Y$ is a Hamiltonian $S^1$-manifold whose
  Hamiltonian function $H$ takes $Y$ onto $I_{\eg}=(\la-\eg,\la+\eg)$
  and $\la$ is the only critical value of $H$.  Moreover we require
  that if $\eg'<\eg$ the restriction $(H^{-1}(I_{\eg'}), H, \eg')$ is
  identified with $(Y,H,\eg')$.  Two cobordisms are equivalent,
  $(Y,H,\eg)\sim(Y',H',\eg')$, if and only if there is
  $\eg''<\min(\eg,\eg')$ such that $(Y,H,\eg'')$ and $(Y',H',\eg'')$
  are isomorphic.  That is, there is an
  isomorphism \[f:H^{-1}(I_{\eg''})\to (H')^{-1}(I_{\eg''}).\] A
  \textbf{critical germ $G(\la,\eps_0)$ at $\la$} is an equivalence
  class in $\hs$ of such tuples.
	
\item[(ii)] Similarly, consider tuples $(Y,H, \eg)$ as above where the
  only critical value of $H$ is its minimum (maximum) value $\la$.
  Thus $H(Y)=[\la,\la+\eg)$ ($H(Y)=(\la-\eg,\la]$).  we have got a similar
  equivalence relation between them as before.  An equivalence class
  $m(\la,\eps) (M(\la,\eps))$ is called a \textbf{minimal (maximal)
    germ} at $\la$.	
\end{itemize}
\end{defn}

We will often refer to critical, maximal or minimal germs just as
germs.  Note that if $\dg<\eps$, there is a natural
restriction map $G(\la,\eps)\to G(\la,\dg)$.
The triples $(Y,H,\eg)$ as in Definition \ref{d:1} (ii) are
neighborhoods of the maximal and minimal sets.  To see this, first
note that there is a unique maximum component $F_{\max}$ of the fixed
point set \cite{AB}.  Then by using an equivariant version of the
Darboux Theorem applied to points in $F_{\max}$, there is a triple of
the form $(Y,H,\eg)$ for $\eg$ small enough.  Moreover, its maximal
germ is determined uniquely by the symplectomorphism type of
$F_{\max}$ and its normal bundle.  A similar remark applies to the
minimum.

\begin{defn}
  Let $I$ be an open interval.  A \textbf{regular slice} is a tuple
  $(Z,K,I, \om)$ where $K:Z\to I$ is a surjective moment map for a
  free $S^1$-action on the symplectic manifold $(Z,\om)$.  We say that
  two slices $(Z,K,I,\om) \sim (Z',K',I,\om')$ are equivalent if they
  are isomorphic.  We denote by $F(I)$ an equivalence class of such
  slices.
\end{defn}

\begin{defn}
  Let $G(\la,\eps)$ be a germ at $\la$ and let $F(I)$ be a class of
  regular slices for $I = (\la',\la)$.  Let $(Y,H,\eg)\in G(\la,\eps)$ 
  and $(Z,K,I)\in F(I)$.  A \textbf{gluing map} $(\phi,\eg):
  (Y,H,\eg)\to (Z,K,I)$ is given by a pair $(\phi,\eg)$ where $0<\eg<
  \eps$ and $\phi$ is an isomorphism
  \[
  H^{-1}(\la-\eg,\la)\stackrel{\phi}\to K^{-1}(\la-\eg,\la). 
  \]
  Two gluing maps $(\phi,\eg), (\phi',\eg')$ are equivalent if there
  are $\eg''<\min(\eg,\eg')$ and isomorphisms \[f:Y\to Y', g:Z\to Z'\]
  such that the following diagram is commutative.
\begin{equation}\label{cd:1}
	\xymatrix{
	H^{-1}(\la-\eg'',\la)\ar[r]^{\phi}  \ar[d]^f  &
        K^{-1}(\la-\eg'',\la) \ar[d]_g \\ 
	{H'}^{-1}(\la-\eg'',\la)\ar[r]^{\phi'} & {K'}^{-1}(\la-\eg'',\la)
	}
\end{equation}

A \textbf{gluing class $\Phi: G(\la,\eps) \to F(I)$} is an equivalence
class of pairs $(\phi,\eg)$.  This maps are well defined by Diagram
\eqref{cd:1}.  Analogously one can define gluing maps for $F(I)$ and
the germ $G(\la,\eps)$ when $I=(\la,\la')$. Note that this definition
also applies if we substitute the germ $G(\la,\eps)$ by a maximal or a
minimal one.
\end{defn}


Recall that we are interested in building symplectic manifolds out of
germs and regular slices.  We now describe the simplest case where we
can do that.  Suppose we have got $G(\la,\eps)$ and $F(\la,\la')$ a
germ and a class of regular slices, we want to see how to obtain a new
manifold via a gluing class $\Phi:G(\la,\eps)\to F(\la,\la')$.  Choose
representatives $(Y,H,\eg)\in G(\la, \eps), (Z,K,I,\om)\in F(I)$ and
$(\phi,\eg)\in \Phi$.  Then, consider the manifold \[Y\cup_{(\phi,
  \eg)} Z\] obtained by gluing $Y$ and $Z$
along the overlap $(\la,\la+\eg)$, that is
\[
Y\sqcup Z/ \sim \text{ where } x\sim y \iff \phi(x)=y
\]
where \[\phi:H^{-1}(\la,\la+\eg)\to K^{-1}(\la,\la+\eg)\] is the
restricted isomorphism on the interval $(\la,\la+\eg)$.  

If $(\phi',\eg'),(Y',H',\eg'),(Z',K',I)$ is another set of choices,
there exist $0<\eg''<\min(\eg,\eg')<\eps$ and isomorphisms $f,g$ as in
the commutative Diagram~\eqref{cd:1}.  Therefore, by restricting the
gluing maps and $f,g$ to the open interval $(\la,\la+\eg'')$ one gets
that $\phi(x)=y$ if and only if $\phi'(f(x))=g(y)$.  Then $(f,g)$
induces an isomorphism
\begin{equation}\label{eq:iso}
Y\cup_{(\phi, \eg)} Z \stackrel{\cong}\to Y' \cup_{(\phi',\eg')} Z'.
\end{equation}
Denote by \[G(\la,\eps)\cup_\Phi F(I)\] the isomorphism class produced
by this gluing.  Note that the Hamiltonian function on this new
manifold is the one defined by $(H,K):Y\sqcup Z \to (\la-\eps,\la')$
after passing to the quotient.  Therefore this is a well defined
operation in $\hs$.

Conversely; the neighborhoods $(H,K)^{-1}(\la-\eg,\la+\eg)\subset
Y\cup_{(\phi, \eg)} Z$ are isomorphic to $H^{-1}(\la-\eg,\la+\eg)$ for
all $\eg$. Similarly $(H,K)^{-1}(\la,\la')\cong
K^{-1}(\la,\la')$. Then we have got the following lemma.

\begin{lemma}\label{l:1}
Suppose we have given a germ of cobordism $G(\la,\eps)$, a regular slice
$F(\la,\la')$ and a gluing class $\Psi: G(\la,\eps)\to F(\la,\la')$. 
Then we can associate a unique isomorphism
class \[G(\la,\eps)\cup_{\Phi} F(\la,\la')\text{ with Hamiltonian
}(H,K)\] in $\hs$. Moreover, the manifolds  $(H,K)^{-1}(I_\eg)$ and
$(H,K)^{-1}(\la,\la')$ represent the classes $G(\la,\eps)$ and
$F(\la,\la')$.  
\end{lemma}


It is clear that we can apply the same idea to define a gluing of the
form $F(\la',\la)\cup_{\Phi} G(\la,\eps)$.   
Similarly, when we have got maximal and minimal germs $M(\la,\eps)$,
$m(\la',\eps)$ one can construct manifolds 
\[
m(\la',\eps')\cup_\Phi F(\la',\la) \text{ and } F(\la',\la)\cup_\Phi
M(\la,\eps).
\] 
It is important to note that  order of the gluing does not matter
provided $\eps$ is small enough. 
To reconstruct $M$, we would like to define this process more generally.  
We want to glue more general data, as we now explain. 
\textbf{A set of local data} $\Ll$ consists of:

\begin{itemize}
\item A collection $\Cc=\{0=\la_0 < \dots < \la_{s +1}\}$ of critical
  levels. 
\item Germs $G(\la_i, \eps_i)$ at $\la_i$ for all $i=1,\dots, s$, 
minimal and maximal germs $m(\la_0,\eps_0)$, $M(\la_{s+1},\eps_{s+1})$
respectively.   
\item For all $j=0,\dots,s$, equivalence  classes of regular slices
  $F(I_j)$ where $I_j=(\la_{j}, \la_{j+1})$ is a maximal open interval
  of  regular values. 
\item Gluing classes  $\Phi^{+}_{i},\Phi^{-}_{i}$ from
  $G(\la_i,\eps_i)$ to  $F(I_{i})$ and $F(I_{i-1})$ respectively for
  all  $i=1,\dots, s$.  
\end{itemize}

As an example consider $(M,\om,H)\in \hs$.  
Its localisations
\[
H^{-1}(\la-\eg,\la+\eg), \text{ and the free sets \ } H^{-1}(\la,\la') 
\]
for $\eg>0$ and $\la,\la'\in\Cc(M)$ define the germs and regular
slices. Then the \textbf{local data associated to $M$}, denoted by
$\Ll_M$, is the collection of all the isomorphism classes of these
sets and the gluing classes on the overlaps. Now we have got the following
theorem.

\begin{thm}\label{t:glu} 
  Given a set of local data $\Ll$, there exists
  a closed Hamiltonian $S^1$-manifold $M_{\Ll}$ such that its
  associated set of local data is $\Ll$. Moreover, this manifold is
  unique up to isomorphism, this is we have got a one-to-one association
  \[ \{ \text{ Sets of local data }\} \to \hs. \] 
\end{thm}

\begin{proof} 

Consider representatives 
\[(Y_{s+1},H_{s+1},\eg_{s+1})\in M(\la_{s+1}, \eps_{s+1}), \ \ \
(Y_0,H_0,\eg_0)\in m(\la_0,\eps_0)\] 
\[((Y_i,H_i,\eg_i)\in G(\la_i, \eps_i) , \eps_{s+1}), \ \ \
(Z_i,K_i,I_i)\in F(I_i)\] and \[(\phi^{\pm}_i,\eg)\in\Phi^{\pm}.\] 
On the disjoint union
\[
\Mm_\Ll = Y_0 \sqcup Z_1 \sqcup Y_1 \dots \sqcup Z_{s+1}\sqcup Y_{s+1} 
\]
define the equivalence relation by
\[
y\sim_\Ll z \iff  \text{for some } j,  (y\in Y_{j}, z\in Z_{j+1} ,
\phi^-_j (y)=z ) \text{ or } (y\in Y_j, z\in Z_{j} , \phi^+_j (y)=z )  
\]

Define $M_\Ll:=\Mm_\Ll/{\sim_\Ll}$ with Hamiltonian $H_\Ll=[H_0,
K_1,\dots, K_{s+1},H_{s+1}]$. Note that $H_\Ll:M_\Ll \to
[\la_0,\la_{s+1}]$ and then $M_{\Ll}$ is a closed symplectic manifold
whose isomorphism class is completely determined by the local data
$\Ll$, since another set of choices would give an isomorphism as in
Equation \eqref{eq:iso}. 
Finally, to see that the association $\Ll \mapsto M_{\Ll}$ is
one-to-one, we proceed as before, by considering the localisations of
$M_\Ll$ at the critical levels with the Hamiltonian $H_\Ll$. One can
see that this recovers $\Ll$. 
\end{proof}

As a last note, we clarify that for the rest of this paper we will
treat isomorphic manifolds as equal, unless  we specify the contrary. 

\section{Determining Local Data.}\label{s:GS} 

\subsection{Basic Notations.}
The background material for this section can be found in
\cite{MS1} and \cite{GS1}.  Let $(M, H,\om)\in \HS$.  Let $t\in
\bb{R}$ be a regular value of $H$. $S^1$ acts freely on the
level set $H^{-1}(t)$. The orbit space $\red{M}{t}:=H^{-1}(t)/S^1$ is
the \textbf{reduced space} of $M$ at the level $t$.  This space is
symplectic with the reduced symplectic form $\overline{\omega}_t$.
The fibration $\pi : H^{-1}(t)\to \red{M}{t}$ is a principal $S^1$
bundle over $\red{M}{t}$.  Denote its total space by $P_{t}$ or just
by $P$, whenever there is no risk of confusion.  Denote its Euler
class by $e(P)\in H^2(\red{M}{t})$.  The reduced symplectic form
$\overline{\omega}_t$ on $\red{M}{t}$ satisfies
\begin{equation}\label{eq:forms}
\pi^*\overline{\omega}_t = i_t^*\om.
\end{equation}


We now introduce the precise definitions of common relations in
symplectic geometry that we use in this paper. Consider two symplectic
forms $\om_{0}, \om_{1}$ on a manifold $X$. These forms are said to be
\textbf{symplectomorphic} if there is a diffeomorphism $f:X\to X$ such
that $f^{*}\om_{1}=\om_{0}$. A \textbf{deformation} between $\om_{0},
\om_{1}$ is a (smooth) family $\{\om_{s}\}$ of symplectic forms that
join them. A deformation is an \textbf{isotopy} if the elements in the
family $\{\om_{s}\}$ all lie in the same cohomology class. It is well
known (Moser's lemma) that two symplectic forms are isotopic if and
only if there is a family of diffeomorphisms $\{h_{s}\}$ on $X$ such
that $h^{*}_{s}\om_{s}=\om_{0}$ and $h_{0}=id$. The concepts of
isotopy and deformation are in general not equivalent (cf. Example
13.20 in \cite{MS1}), but for some special cases as we will see in
Theorem \ref{t:dusa}, they agree. For the objectives of the present
work, manifolds where these two properties agree will be a key
ingredient as we will see in Lemma \ref{c:equiv}.

When the manifolds are equipped with circle actions, we will make the
natural assumption that all the deformations and maps are
$S^{1}$-equivariant.

\subsection{The isomorphism type of regular slices.}\label{ss:free} As
before, let $(M,H,\om)\in \HS$. Suppose $\la'<\la$ are two consecutive
critical values of $H$. The interval $J=(\la',\la)$ is a maximal set
of regular values of $H$. $S^{1}$ acts freely on the open set
$H^{-1}(J)$. Let $I=[t_0,t_1]\subset J$ and let $(\red{M}{t},
\overline{\om}_{t})$ denote the reduced spaces for $t\in I$. The
classic Duistermaat-Heckman theorem establishes that the reduced
spaces $\red{M}{t}$ are all identified with $B:=\red{M}{t_0}$ for a
fixed value $t_0\in I$. The diffeomorphism type of $H^{-1}(I)$ is then
given by $P_{t_0}\times I$. The forms $\overline{\omega}_t$ define a
path of symplectic structures on $B$, whose cohomology classes satisfy
\begin{equation}\label{eq:DH}
 [\overline{\omega}_t]=[\overline{\omega}_{t_0}]
  +(t-t_0)e(P), t\in I. 
\end{equation}
The symplectic structure of $H^{-1}(I)$ is completely characterized as
follows.
\begin{lemma}[Proposition 5.8 in \cite{MS1}]\label{l:DH}  
  The symplectic manifold $(H^{-1}(I), \omega)$ is determined by the
  bundle $P\to B$ and the family $\{\overline{\omega}_{t\in I}\}$ up
  to an equivariant symplectomorphism.
\end{lemma}

Our aim in this section is to understand the isomorphism class of
$H^{-1}(I)$ in terms of less information. Suppose that we have got two
invariant symplectic forms $\om$ and $\om'$ on $H^{-1}(I)$. Equation
\eqref{eq:DH} shows that the paths defined by the cohomology classes
of the reduced forms satisfy
$[\overline{\omega}_t]=[\overline{\omega}'_t]$ in $B$ provided
$[\overline{\omega}'_{t_0}] = [\overline{\omega}_{t_0}]$. This,
motivates the following.

\begin{definition}
  Let $B$ be a manifold. We say that two paths $\overline{\omega}_t,
  \overline{\omega}'_t, t\in I$ of symplectic forms are
  \textbf{weakly-equivalent} if $[\overline{\omega}_t] =
  [\overline{\omega}'_t]$ for all $t\in I$. This is equivalent to
  require $[\overline{\omega}'_{t_0}]=[\overline{\omega}_{t_0}]$ and
  impose Equation \eqref{eq:DH} for some class $e(P)$. We denote by
  $\breve{\overline{\omega}}_t$ the weak equivalence class of the path
  $\{\overline{\omega}_t\}_{t\in I}$.
\end{definition}

In particular, any two invariant symplectic forms $\omega, \omega'$ in
$H^{-1}(I)$ yield weakly-equivalent paths of reduced forms if their
symplectic reductions at $t_0$ agree.  What we want to know is when
any two families of forms $[\overline{\omega}_t],
[\overline{\omega}'_t]$ in $B$ satisfying Equation \eqref{eq:DH} lift
to isomorphic forms on $H^{-1}(I)$.  Equivalently we want to
determine the isomorphism type of $P\times I$ in therms of these paths
of forms. Suppose that there is an equivariant diffeomorphism
$\phi:P\times I\to P\times I$ such that $\phi^*(\om)=\om'$.  $\phi$
defines a family $\phi_t:P_t\to P_t$ by restricting to the regular
level at $t$. These maps descend to the quotient $\red{\phi}{t}:B\to
B$, satisfying the commutative diagram
\begin{equation}\label{eq:diagram1}
  \xymatrix{
    P_t \ar[d]_\pi \ar[r]^{\phi_t} & P_t \ar[d]_\pi\\
    B \ar[r]^{\red{\phi}{t}} & B
  }
\end{equation}
The maps $\red{\phi}{t}$ satisfy $(\red{\phi}{t})^*
(\overline{\omega}_t) = \overline{\omega}'_t$ and preserve the Euler class
$e(P)$ for all $t\in I$.

Conversely, given any smooth family of maps $\red{f}{t}
: B \to B$ such that $(\red{f}{t})^*\overline{\omega}_t=
\overline{\omega}'_t$ and $(\red{f}{t})^* e(P)=e(P)$, it lifts to a
family of isomorphisms $f_t:P \to P$ that make Diagram
\eqref{eq:diagram1} commute. Therefore, these maps bundle together to
a diffeomorphism $f:P\times I\to P\times I$ such that $f^*(\om)=\om'$.

For simplicity of notation, assume that $I=[0,1]$. 

\begin{definition}\label{d:rel} 
  Two weakly-equivalent families $\{ \overline{\omega}_{t} \}_{t\in
    I}$ and $\{\overline{\omega}'_{t}\}_{t\in I}$ of symplectic forms
  on $B$ are said to be \textbf{equivalent} if there is a smooth
  family $\{\overline{\omega}_{s,t}\}$ of symplectic forms such that
  \begin{equation}\label{eq:cond}
    \frac{d}{ds}[\overline{\omega}_{s,t}]=0 \  
    \text{ and }  \overline{\omega}_{0,t}=\overline{\omega}_t,
    \overline{\omega}_{1,t}=\overline{\omega}'_t, 
  \end{equation}
  for $0\le t,s \le 1$. 
\end{definition}

If $\{\overline{\omega}_{t}\}$ is equivalent to
$\{\overline{\omega}'_{t}\}$, then there is an isotopy $\Omega_s$ of
forms on $P\times I$ that lifts $\overline{\omega}_{s,t}$ and such
that $\Omega_0=\om, \Omega_1=\om'$.  By Moser's lemma, we have got a
family of maps $f_s:P\times I \to P\times I$ such that
$f_s^*\Omega_s=\om$. Therefore the isomorphism type of $P\times I$
depends exclusively on the equivalence class of
$\{\overline{\omega}_t\}$.  Using rigidity (Definition \ref{d:rigid})
we now understand the relation amongst weakly-equivalent and
equivalent.

\begin{lemma}\label{c:equiv}
  Let $(B,\{\overline{\omega}_t\}_{t\in I})$ be a rigid pair. Then,
  any family $\{\overline{\omega}_{t}'\}$ weakly-equivalent to
  $\{\overline{\omega}_{t}\}$ and with $\overline{\omega}_{0}=
  \overline{\omega}'_{0}$, is also equivalent to
  $\{\overline{\omega}_{t}\}$.
\end{lemma} 

\begin{proof}
  We want to see that there is a family of
  symplectic forms $\{\overline{\omega}_{s,t}\}_{0\le s,t, \le1}$ such
  that for each fixed $t$ the path $\overline{\omega}_{s,t}$ is an
  isotopy from $\overline{\omega}_t$ to $\overline{\omega}'_t$.
  First, since the families are weakly-equivalent,
  $\overline{\omega}_{t}$ is cohomologous to $\overline{\omega}_{t}'$
  for all $t$, then for $s\in[0,1]$ the cohomology class of the form
  $s \overline{\omega}'_t + (1-s) \overline{\omega}_t$ is constant
  with respect to $s$. Since $\overline{\omega}_0 =
  \overline{\omega}'_0$, Moser's argument shows that there is an
  $\eg>0$ and a family $\overline{\omega}_{s,t}$ satisfying Equation
  \eqref{eq:cond} for $0\le t \le \eg$. (Compare with Example 3.20 in
  \cite{MS1}).  We want to see that we can take
  $\eg=1$.  To see this, define
  \[
  \Dd=\{ T : \exists \ \{\overline{\omega}_{s,t}\}_{0\le s\le 1,
    0\le t\le T} \text{ satisfying Equation \eqref{eq:cond}} \}\subset
  [0,1]. 
  \]
  We claim that $\Dd$ is $[0,1]$.  We will do this by proving that
  $\Dd$ is open and closed. Let $T\in \Dd$. Since
  $\overline{\omega}_T$ is isotopic to $\overline{\omega}'_T$, we can
  assume that $\overline{\omega}_T=\overline{\omega}_T'$. Then, by the
  same argument as when $T=0$, there is an $\eps>0$ such that
  $T+\eps\in \Dd$. Thus $\Dd$ is open.

  To see that $\Dd$ is closed, take $T$ be such that $0<T\le 1$, and
  $T-\eg\in \Dd$ for every $\eg>0$ small. The path
 \begin{equation} 
   \al_s=
   \begin{cases} \overline{\omega}_{v(s)} & v(s)=(1-2s)T, 0\le s\le\half \\
     \overline{\omega}'_{u(s)} & u(s)=(2s-1)T, \half\le s \le 1
   \end{cases}
 \end{equation}
 is a deformation between $\overline{\omega}_T$ and
 $\overline{\omega}'_T$. Since $(B,\{\overline{\omega}_t\})$ is rigid
 and $[\overline{\omega}_T]=[\overline{\omega}'_T]$, $\ag_s$ can be
 homotoped through deformations with fixed endpoints to an
 isotopy. Let $\bg_s$ such isotopy with $\bg_0=\overline{\omega}_T$
 and $\bg_1=\overline{\omega}'_T$. Again, by Moser's argument, we can
 extend $\bg_s$ in a small neighborhood of $T$.  That is, for an
 $\eg>0$ there is a family $\bg_{s,T-\eg}$ of symplectic forms which
 is homotopic to $\bg_s$.

 On the other hand, because of the hypothesis on $T$, we have
 $T-\eg\in \Dd$.  Thus there is a family $\overline{\omega}_{s,t}$
 that satisfies Equation \eqref{eq:cond}, for all $t\in[0,T-\eg]$.
 The concatenation of the two isotopies $\overline{\omega}_{s,T-\eg}$,
 $-\bg_{s,T-\eg}$ defined by
 \[
 \ga_{s} :=
 \begin{cases}
   \overline{\omega}_{s,T-\eg} & 0\le s\le 1\\
   \bg_{2-s,T-\eg} & 1\le s \le 2
 \end{cases}
 \] 
 (after smoothing) defines a loop at $\overline{\omega}_{T-\eg}$ in
 the space of symplectic structures $\Ss(a)$ with fixed cohomology
 class $a=[\overline{\omega}_{T-\eg}]$.
The fibration 
\[
\Symp(B,\overline{\omega}_{T-\eg})\cap \Diff_0(B) \to
\Diff_0(B)\stackrel{\pi}\to \Ss(a), \ \ \ \pi : f\mapsto
f^*(\overline{\omega}_{T-\eg}),
\]
gives a lift $\{f_{s}\}_{s\in[0,2]}$ in $ \Diff_0(B)$ such that
$f_0=id$ and $f_{s}^*(\ga_s)=\overline{\omega}_{T-\eg}$ for all $s$.
By rigidity hypothesis, the fibre
$\Symp(B,\overline{\omega}_{T-\eg})\cap \Diff_0(B)$ is path connected,
we can assume that $f_2=id$ as well.

The map \[h_s:=f_{s}\circ (f_{2-s})^{-1}, \ 0\le s\le 1\] is such
that \[h_0=h_1=id, \ \ h_s^*(\ga_{2-s})=\ga_s\] that is,  
\[
h_s^*(\bg_{2-s,T-\eg}) = \overline{\omega}_{s,T-\eg}.
\]
Therefore, the new family
\[ \hat{\bg}_{s,t} := h_s^*(\bg_{s,t})\] satisfies Equation
\eqref{eq:cond} for $0\le s\le1, T-\eg\le t \le T$ and agrees with
$\overline{\omega}_{s,t}$ at $t=T-\eg$. After smoothing, we see that
$\overline{\omega}_{s,t}$ can be extended to all $t\le T$ via
$\hat{\bg}_{s,t}$. Then $T\in \Dd$. This proves that $\Dd$ is closed.
 
We would like to emphasise where our argument fails if
$\Symp(B,\overline{\omega}_{T-\eg})\cap \Diff_0(B)$ is not path
connected. In this case, we cannot consider $f_2$ to be the identity
and then $h_0\ne id$. For our boundary conditions, we need the
extension $\hat{\bg}_{0,t}$ to agree with $\overline{\omega}_t$ for
$T-\eg\le t\le T$. We cannot conclude that if $h_0$ is not the
identity. 
\end{proof}

Combining all the results above, yield the most important result of this section.

\begin{theorem}\label{p:equiv}
  Let $(M,H,\omega)$ be a Hamiltonian $S^1$-manifold with proper
  Hamiltonian and let $\lambda',\lambda$ be any two consecutive
  critical values. Let $I=[t_0,t_1]\subset (\lambda',\lambda)$, and
  suppose that under the identification $\red{M}{t}=\red{M}{t_0}$,
  $(\red{M}{t_0}, \{\overline{\omega}_t\}_{t\in I})$ is rigid. Then
  any other invariant form $\omega'$ on $H^{-1}(I)$ is isomorphic to
  $\omega$ provided its reduction at $t_0$ is (diffeomorphic to)
  $\overline{\omega}_{t_0}$. In other words, the isomorphism class of
  the regular slice $H^{-1}(I)$ is determined up to isomorphism by the
  bundle $P\to \red{M}{t_0}$, the family of cohomology classes
  $[\overline{\omega}_t]_{t\in I}$ and the initial form
  $\overline{\omega}_{t_0}$.
\end{theorem}

\begin{proof}
  Since any two invariant forms on $H^{-1}(I)$ with
  $\overline{\omega}'_{t_0} =\overline{\omega}_{t_0}$ yield
  weakly-equivalent paths of forms the family of reduced forms
  $\{\overline{\omega}'_t\}_{t\in I}$ is weakly-equivalent to
  $\{\overline{\omega}_t\}_{t\in I}$, Lemma \ref{c:equiv} shows they
  are also equivalent, therefore $\omega, \omega'$ are isomorphic.
\end{proof}

%
%

\subsection{Germs at critical levels.}\label{ss:critical}

We now describe the germ $G(\la, \eps)$ when $H^{-1}(\la)$ contains
only fixed points of index 2 and the action is semi-free. We will also
discuss the change in the fixed point data.

\subsubsection{Smooth and symplectic structure on the critical reduced
  space.} \label{sss:smooth} Let $H^{-1}(\lambda)$ be a level with
only index 2 fixed-point components for $(M,H,\omega)$. In general
$H^{-1}(\lambda)$ is a singular space and thus $\red{M}{\la}$ is not
naturally smooth. Although it is well known that we can provide this
space with a smooth symplectic manifold structure \cite{GS1}, one can
put a canonical, independent of choices, smooth structure on
$\red{M}{\la}$ using \emph{grommets} near the fixed point set, as in
Tolman-Karshon \cite{KaTo-Ce}. In what follows we use a
re-interpretation due to McDuff \cite[\S 3.2]{Mc-So09}.  Let $F\subset
H^{-1}(\lambda)$ denote the fixed point set, and for simplicity assume
it is connected of codimension $2k$. The equivariant Darboux-Weinstein
theorem states that germs of neighborhoods $U$ in $M$ containing $F$
are isomorphic. One can identify any such $U$ with the normal bundle
$N_{F/M}$ of $F$ in $M$. The isomorphism class of $U$ is determined by
the restriction of the symplectic structure to $F$ and the isomorphism
class of the normal bundle \cite[3.30]{MS1}. To give a symplectic
chart to the reduced space, consider Darboux coordinates of $F$ near
$x\in F$, $(w_1, \dots w_{n-k})$ in $\mathbb{C}^{n-k}$, and let
$(z_1,\dots, z_k)$ denote Darboux coordinates in the normal directions
of $U$, so that $H$ is given in the standard form
$\lambda-|z_1|^2+|z_2|+\dots +|z_k|^2$. This local representation is
unique up to equivariant isotopy. The map (grommet)
\begin{equation}\label{eq:smooth} 
\ell:\mathbb{C}^{n-1}\to
  H^{-1}(\lambda); \ \ \ 
 (w_1, \dots, w_{n-k},z_2,\dots,z_k)\mapsto
  (w_1, \dots, w_{n-k},\left(\sum|z_j|^2\right)^{\half}, z_2,\dots,
  z_k) 
\end{equation}
intersects each orbit inside the singular set $H^{-1}(\lambda)$
exactly once. This defines a homeomorphism from the reduced space
$\red{M}{\lambda}$ to $\mathbb{C}^{n-1}$, which gives a chart near the
image $\overline{F}$ of $F$ in $\red{M}{\lambda}$. The non-singular
locus $\red{M}{\lambda} - \overline{F}$ is naturally smooth. Patching
these pieces together, we get a natural smooth structure independent
of choices. The symplectic structure $\omega$ on $M$ descends to
$\red{M}{\lambda} - \overline{F}$, which together with the symplectic
form using \eqref{eq:smooth} near $\overline{F}$ gives a symplectic
structure $\red{\omega}{\lambda}$ on $\red{M}{\lambda}$ which depends
exclusively on $(\red{M}{\lambda},\overline{F},\omega)$.  It follows
from this construction that all the reduced spaces $\red{M}{t}$ for
$t\in(\lambda-\eps,\lambda)$ are diffeomorphic to $\red{M}{\lambda}$
for suitable small $\eps$. 

\subsubsection{Cobordism and the change of fixed point
  data.}\label{sss:cobord} 
Let $(M,\og)$ be a symplectic manifold and let $S\subset M$ be a
closed submanifold. For $\eg>0$ small enough, we denote by $\Bl{M}{S}$
and $\be:\Bl{N}{S}\to N$ the $\eg$-symplectic blow up of $M$ along $S$
and the blow-down map as defined by Guillemin and Sternberg
\cite{GS1}. The construction of the manifold $\Bl{M}{S}$ depend on
several choices\footnote{To be precise, an embedding of a ball bundle
  and a connection on a principal $U(n)$-bundle.}, but its
diffeomorphism type is independent of them.  $\Bl{M}{S}$ admits a
blow-up symplectic structure denoted by $\widehat{\om}(\eg)$.  This
form is not independent of the choices\footnote{Not even when one
  blows-up a point.  This is because there is a non-compact family of
  choices.}, but its \textsl{germ} of isotopy classes is. That is, if
$\widehat{\om}'(\eg)$ is a form obtained by making different choices,
then for some $\eps_0$ small enough there exist a smooth family
$f_\eg\in \Diff(M)$ such that $f_\eg^* (\widehat{\om}(\eg))=
\widehat{\om}'(\eg)$ for all $0<\eg<\eps_0$. With this in mind, assume
we have got the following data.
\begin{enumerate}
  \item A compact symplectic manifold $(\overline{M},\overline{\om})$.
  \item A symplectic submanifold $\overline{F}\subset \overline{M}$.
  \item A principal $S^{1}$-bundle $\pi:P\to \overline{M}$.
  \item A connection 1-form $\ag$ on $P$.
\end{enumerate}

Then it is possible to create a cobordism $(Y({\la}),H,\eg)\in \HS$
having the following properties. 

\begin{enumerate}
\item $H$ maps $Y({\la})$ onto $I=(\la-\eg,\la+\eg)$.
\item For all $t>0$, $H^{-1}(\la- t)$ is equivariantly
  diffeomorphic to $P$. 
\item For all $t>0$ the symplectic reduction $Y({\la})$ at the level
  $\la -t$ is symplectomorphic to $\overline{M}$ with symplectic form
  $\overline{\og}-td\ag$.\footnote{$d\al$ is a 2-form on $P$, but it
    descends to $\overline{M}$. This is the form that we consider
    here.}
\item $\la$ is a critical level of $H$ of index 2.
\item For all $t>0$ the reduction of $Y(\la)$ at $\la+t$ is the blow
  up $\Bl{\overline{M}}{F}$ of $\overline{M}$ along $F$ with symplectic
  structure
  \begin{equation}\label{eq:localform}
    \widehat{\om}(t) +\bg^*(td\ag).
  \end{equation}
  Here $\widehat{\om}(t)$ is the $t$-blow up form and
  $\bg:\Bl{\overline{M}}{F}\to \overline{M}$ is the blow down map.
  \item The fixed point set at  $\la$ is $F$. 
\end{enumerate}



\begin{theorem}[Guillemin-Sternberg]\label{t:bordism}
  Let $(M,K,\om)\in \HS$ and $\la$ be an index 2 critical value. Let
  $I=(\la-\eps, \la+\eps)$ be a sufficiently small interval in
  $\bb{R}$.  Then, the open submanifold $K^{-1}(I)$ is equivariantly
  symplectomorphic to the manifold $(Y(\la),H,\eps)$. Moreover, the
  germ of diffeomorphism of $M$ near $\la$ only depends on the fixed
  point data $(\red{M}{\la},\red{F}{},\overline{\omega}_\la)$,
  $e(P_{\la})\in H^2(\red{M}{\la})$.
\end{theorem}

Theorem \ref{t:bordism} describes the change in the fixed point data
after crossing a critical level $\la$ with index 2 fixed point
sets. A similar analysis applies for coindex 2. We will include in
Lemma \ref{l:iso-germ} details of the proof of the uniqueness
statement, since we will need it later.

Let $N_{F/M}$ be the normal bundle of the fixed point submanifold $F$
in $M$. This bundle decomposes as
\begin{equation}
N_{F/M}=N^+_{F/M}\oplus N_{F/M}^-\label{eq:isobundle0}
\end{equation}
taking positive and negative directions. A small neighborhood of the
zero section of $N_{F/M}$ is $S^1$-isomorphic to a neighborhood
$U\cong U^+\oplus U^-$ of $F$ in $M$. Denote by $\red{X}{t}$ the
symplectic reduction of $U^{-}$ at $t$ if $t\le\la$ or the reduction
of $U^{+}$ at $t$ if $t>0$. The diffeotype of the triple
$(P_t,\red{M}{t}, \red{X}{t})$ depends smoothly on $t$ for
$t\in(\la-\eg,\la]$ and it is constant. Denote it by $(P_-,\red{M}{-},
\red{X}{-})$ where $\red{X}{-}=\red{F}{\la}$.  Similarly, we denote
the common diffeotype for $t\in (\la,\la+\eg)$ by $(P_+,\red{M}{+},
\red{X}{+})$.

We now remark an important relation of the bundles
\cite[pp. 516]{GS1}. By analysing the negative directions, there is an
isomorphisms of $S^1$-bundles
\begin{equation}\label{eq:isobundle1}
  P_\lambda|_{\red{F}{}}\cong \mathrm{Frame}(N_{F/M}^{-})
\end{equation}
where $P_\lambda|_{\red{F}{}}$ is the restriction of $P_\lambda$ to
the fixed point set in the reduced space, and
$\mathrm{Frame}(N_{F/M}^{-})$ is the frame bundle of
$N_{F/M}^{-}$. Similarly there is an isomorphism of bundles
\begin{equation}
  N_{\red{F}{}/\red{M}{\lambda}}\cong N_{F/M}^+\label{eq:isobundle2}.
\end{equation}
Equations \eqref{eq:isobundle0},\eqref{eq:isobundle1} and
\eqref{eq:isobundle2} show that the normal bundle $N_{F/M}$ is
determined by $N_{\red{F}{}/\red{M}{\lambda}},\overline{F}$ and
$P_\lambda$. 

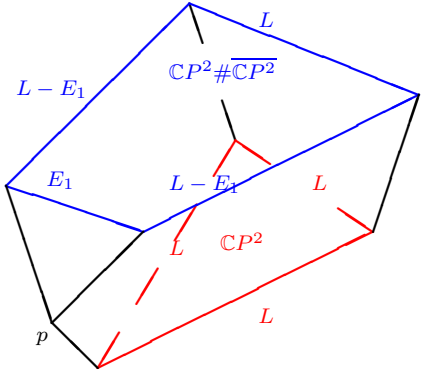
\begin{figure}
\begin{center}
\setlength{\unitlength}{2500sp}%
\begingroup\makeatletter\ifx\SetFigFont\undefined%
\gdef\SetFigFont#1#2#3#4#5{%
 \reset@font\fontsize{#1}{#2pt}%
  \fontfamily{#3}\fontseries{#4}\fontshape{#5}%
  \selectfont}%
\fi\endgroup%
\begin{picture}(4000,4000)(0,-3500)
\put(2500,-800){\makebox(0,0)[lb]{\smash{{\SetFigFont{8}{8}{\rmdefault}
{\mddefault}{\updefault}{\color[rgb]{0,0,1}$\CP{2}\#\overline{\CP{2}}$}%
}}}}
\thicklines
{\color[rgb]{0,0,0}\multiput(2701,-61)(319.96462,-959.89387){2}{\line( 1,-3){130.035}}
}%
{\color[rgb]{1,0,0}\put(1801,-3661){\line( 2, 1){2700}}
\multiput(3151,-1411)(-379.47841,-632.46402){4}{\line(-3,-5){211.565}}
\multiput(3151,-1411)(1007.85429,-671.90286){2}{\line( 3,-2){342.146}}
}%

{\color[rgb]{0,0,0}\put(1351,-3211){\line( 1,-1){450}}
}%
{\color[rgb]{0,0,1}\put(2251,-2311){\line( 2,
    1){2700}}\put(901,-1861){\line( 1, 1){1800}} 
\put(2701,-61){\line( 5,-2){2250}}\put(901,-1861){\line( 3,-1){1350}}
}%
{\color[rgb]{0,0,0}\put(4501,-2311){\line( 1, 3){450}}
}%
{\color[rgb]{0,0,0}\put(1351,-3211){\line( 1, 1){900}}
}%
{\color[rgb]{0,0,0}\put(901,-1861){\line( 1,-3){450}}
}%
\put(1200,-3400){\makebox(0,0)[lb]{\smash{{\SetFigFont{8}{8}{\rmdefault}{\mddefault}{\updefault}{\color[rgb]{0,0,0}$p$}   
}}}}
\put(1300,-1861){\makebox(0,0)[lb]{\smash{{\SetFigFont{8}{8}{\rmdefault}{\mddefault}{\updefault}{\color[rgb]{0,0,1}$E_{1}$}%
}}}}
\put(3376,-286){\makebox(0,0)[lb]{\smash{{\SetFigFont{8}{8}{\rmdefault}{\mddefault}{\updefault}{\color[rgb]{0,0,1}$L$}%
}}}}
\put(1000,-961){\makebox(0,0)[lb]{\smash{{\SetFigFont{8}{8}{\rmdefault}{\mddefault}{\updefault}{\color[rgb]{0,0,1}$L-E_{1}$}%
}}}}
\put(2500,-1900){\makebox(0,0)[lb]{\smash{{\SetFigFont{8}{8}{\rmdefault}{\mddefault}{\updefault}{\color[rgb]{0,0,1}$L-E_{1}$}%
}}}}
\put(3900,-1900){\makebox(0,0)[lb]{\smash{{\SetFigFont{8}{8}{\rmdefault}{\mddefault}{\updefault}{\color[rgb]{1,0,0}$L$}%
}}}}
\put(3376,-3211){\makebox(0,0)[lb]{\smash{{\SetFigFont{8}{8}{\rmdefault}{\mddefault}{\updefault}{\color[rgb]{1,0,0}$L$}%
}}}}
\put(3000,-2500){\makebox(0,0)[lb]{\smash{{\SetFigFont{8}{8}{\rmdefault}{\mddefault}{\updefault}{\color[rgb]{1,0,0}$\CP{2}$}%
}}}}
\put(2500,-2550){\makebox(0,0)[lb]{\smash{{\SetFigFont{8}{8}{\rmdefault}{\mddefault}{\updefault}{\color[rgb]{1,0,0}$L$}%
}}}}
\end{picture}
\end{center}
\caption{The cobordism around an isolated fixed point $p$ of index
  2. Here we assume that the minimum is also isolated. The base of the
  figure represents  $\CP{2}$ and the top the blow-up
  $\CP{2}\#\overline{\CP{2}}$. If the fixed point is of index 4  this
  cobordism is up-side down.} 
\label{f:1}
\end{figure}
%

Now we address the relation amongst the principal bundles $P_+$ and
$P_-$. This is quite subtle.  By Theorem \ref{t:bordism},
$\red{M}{+}$ is the blow up of $\red{M}{-}$ along $\red{X}{-}$.  Let
$\bg:\red{M}{+} \to \red{M}{-}$ be the blow down map.  $\bg$ restricts
to a diffeomorphism $\red{M}{+}-\red{X}{+}$ onto
$\red{M}{-}-\red{X}{-}$, and when restricted to $\red{X}{+}$ it is a
fibration
\[ \bg: \red{X}{+} \to \red{X}{-}\] whose fibres are all diffeomorphic
to $\CP{k-1}$. Here $2k$ is the codimension of $\red{X}{-}$ in
$\red{M}{-}$.
Denote by $L'$ the line bundle on $\red{M}{+}$ whose Chern class is
dual to its codimension 2 submanifold  $\red{X}{+}$. Let $L$ be the
circle bundle associated to $L'$.  
Then we have \cite[Formula 13.3]{GS1}
\begin{equation}\label{eq:euler}
P_{+}= \bg^{*}(P_{-})\otimes L 
\end{equation}
as circle bundles over $\red{M}{+}$. Since $L$ is trivial on
$\red{M}{+}-\red{X}{+}$, then $P_{-}\cong P_{+}$ on
$\red{M}{+}-\red{X}{+}$. The construction of $L$ depends on the normal
bundle of $\red{X}{+}$ in $\red{M}{+}$, and hence on the pair
$(\red{M}{-},\red{X}{-})=(\red{M}{\la},\red{F}{\la})$. Then one is
tempted to describe $P_+$ in terms of $P_-,(\red{M}{-},\red{F}{\la})$
. We now see examples on which the relation between $P_{-}$ and
$P_{+}$ is easy to depict.

\begin{example}\label{ex:1} \rm Suppose $F$ is an isolated fixed
  point $p\in M$ of index 2.  The relation between the Euler
  classes $e(P_-)$ and $e(P_+)$ is clear. Let $\beta: \overline{M}_+\to
  \overline{M}_-$ denote the blow down map that collapses  the
  exceptional divisor $E$ to $p\in \overline{M}_-$. Then we have
  \[
  e(P_+)=\bg^{*}e(P_-) + \PD(E).
  \]
  In the case $p$ has coindex $2$, the relation is inverted 
  \[
  \bg^{*}e(P_+)=e(P_-) + \PD(E),
  \]
  since now $\overline{M}_-$ blows down to $\overline{M}_+$.
\end{example}

\begin{example}\label{ex:2} 
  Suppose $\mathrm{codim } \ F=4$ and that the index of $F=2$.  By
  Lemma 5 in McDuff \cite{Mc2} all the reduced spaces $\red{M}{t}$ for
  $t\in I, t\ne\la_{1}$ are diffeomorphic, to $B_-$ say.  In
  particular $B_+$ is diffeomorphic to $B_-$.
  Moreover \[e(P_+)=e(P_-)+\PD(F)\] where $F$ is embedded in $B_-$ as
  before. If $\dim M=6$ this case applies when $F$ is a surface.
\end{example}

These two examples show, that in special circumstances the bundles $P$
are determined by the rest of the fixed point data at and below the
critical level $\lambda$ and thus they might be discarded as input
information for the classification problem.

\subsubsection{The critical germ at non-extremal critical levels.} It is pointed
out in \cite[pp. 514]{GS1} that the symplectomorphism type of the
germs is independent of choices and it depends exclusively in the
fixed point data. We include the argument in \cite{Mc-So09} to show that
any two germs with the same fixed point data at a critical value are
isomorphic. We then adapt all these ideas to show, that using rigidity
we can remove more information from the fixed point data without
altering the isomorphism type of the germ.  Before we state the
results, we recall that we are using the smooth structure for reduced
spaces as in \S\ref{sss:smooth} and that to avoid complications we are
assuming each critical level to be simple. Please refer to \S
\ref{ss:non-simple} for the non-simple case.

\begin{definition}\label{d:sfpd1}
  Let $(M_i,H_i,\omega_i), i=1,2$ denote two Hamiltonian
  $S^1$-manifolds.  We say that they have the \textbf{same fixed point
    data} at a \textbf{non-extremal critical value} $\lambda$ if
  \begin{enumerate} 
  \item $\lambda$ is a non-extremal critical value for both $H_i$.  
  \item\label{d:symp} There is a symplectomorphism of reduced spaces
    $\phi:(\red{M}{1,\lambda}, \overline{F}_{1,\lambda},
    \overline{\omega}_{1,\lambda})\to(\red{M}{2,\lambda}, \overline{F}_{2,\lambda},
    \overline{\omega}_{2,\lambda})$ such that $\phi^*e(P_{2,\lambda})=e(P_{1,\lambda})$.
  \item The index functions are related by $i_{2,\lambda}\circ \phi =i_{1,\lambda}$.
  \end{enumerate}
  We say that they have the \textbf{same fixed point data at the
    minimum} if $\lambda$ is the minimum value for both $H_i$ and
  there is a symplectomorphism $\phi:(F_{1,\min},
  \omega_{F_{1,\min}})\to (F_{2,\min},\omega_{F_{2,\min}})$ such that
  $\phi^* (N_{F_{2,\min}/ M_2})=N_{F_{1,\min}/M_1}$. Similarly for the
  maximum.
\end{definition}

The following two lemmas deal with the isomorphism type of the
germs. The first one is a minor generalisation of the result in
\cite{Mc-So09}. The second one is just a small adaptation of the proof
of Theorem 13.1 in \cite{GS1}, where it is shown that two particular
manifolds with the same fixed point data are isomorphic.

\begin{lemma}\label{l:iso-germ}
  Let $(M_i,H_i,\omega_i), i=1,2$ denote two (possibly open)
  Hamiltonian $S^1$-manifolds with proper Hamiltonians $H_i$. Suppose
  that $\lambda$ is a common critical value of $(M_i,H_i, \omega_i
  ), i=1,2$ with only fixed point components of index 2, and they have
  the same fixed point data at $\lambda$.  Then, there is an
  $\eps_0 > 0$ such that $H_1^{-1}(\lambda-\eps, \lambda + \eps)$ is
  isomorphic to $H_2^{-1}(\lambda - \eps, \lambda + \eps)$, for all
  $\eps \leq \eps_0$. 
\end{lemma}

\begin{proof}
  Assume for now that the fixed point set is a single isolated fixed
  point $p_i$. Let $\chi:U_i\to U\subset\mathbb{C}^n$ denote a Darboux
  chart near $p_i$ so that $H_i\circ\chi^{-1}$ are in standard form
  $\lambda-|z_1|^2 + \sum_{1<j\leq n}|z_j|^2$. Since the fixed point data is
  the same near $\lambda$ there is a symplectomorphism
  $\phi:(\red{M}{1,\lambda},
  \overline{\omega}_{1,\lambda})\to(\red{M}{2,\lambda},
  \overline{\omega}_{2,\lambda})$. In the local model we can assume
  that this map can be isotoped to the identity, by shrinking the
  neighborhoods $U_i$ is necessary. Thus using the charts in the
  reduced spaces given by \eqref{eq:smooth} we can assume that
  \(\ell^{-1}\circ\chi_2\circ\phi\circ\chi_1^{-1}\circ\ell\) is the
  identity near the fixed point set.
  
  Claim 1. For $\eps>0$ small enough, the two manifolds
  $H_i^{-1}(\lambda-\eps,\lambda+\eps)\backslash U_i$ are
  equivariantly diffeomorphic.\footnote{Compare this and what follows
    with proof of Theorem 13.1 in \cite{GS1}, they denote this piece
    as $Q$, note that their hypothesis and statement is slightly
    different.}
  
  This claim follows essentially from the fact that on the complement
  $H_i^{-1}(\lambda-\eps,\lambda+\eps)\backslash U_i$ the action is
  free, and since both manifolds have the same symplectic reduction at
  zero, it follows from Duistermaat-Heckman theorem.  To make this
  explicit, choose an $S^1$-invariant and compatible almost complex
  structure on $M_i$ so that near $p$ they agree with $\chi^*(J_0)$,
  where $J_0$ the standard complex structure in $\mathbb{C}^n$. Let
  $g_i$ denote the metric induced by $\omega_i, J_i$ and let
  $\tau_{i,t}$ be the downwards gradient flow with respect to the
  metrics $g_i$. By taking $\eps>0$ small enough we can assume that
  the neighborhoods $U_i$ are compatible with the gradient flow. This
  is we can assume that the set $U_i\cap
  H^{-1}_i(\lambda-\eps,\lambda+\eps)$ is not empty and that it
  contains all the orbits in $H^{-1}_i[\lambda,\lambda+\eps)$ whose
  downward gradient flow converges to the fixed points. By
  re-parametrising if needed one can assume that $\tau_{i,t}$ takes
  any level set $H_i^{-1}(c)$ to $H_i^{-1}(c-t)$. These flows are
  non-singular away from $p_i$, so that outside $U_1$ one can define
  the equivariant diffeomorphism $\tau(x,t) = \tau_{2,t}\circ
  \hat{\phi} \circ \tau^{-1}_{1,t}(x)$, where $\hat{\phi}$ is the
  symplectomorphism of $H_1^{-1}(\lambda) \backslash U_1$ to
  $H_2^{-1}(\lambda) \backslash U_2$ lifting $\phi$. This proves the
  claim.

  Now define $\tau$ on $U_1$ to be the isomorphism $\chi_2^{-1} \circ
  \chi_1$. Then $\tau^*\omega_2=\omega_1$ on $U_1$ by construction on
  $U_1$, but it might not be outside $U_1$.  The family of forms
  $(\tau^*\omega_2)_t, (\omega_1)_t$ induced in the reduced spaces
  $(H_1^{-1}(t) \backslash U_1)/S^1$ define a path of non-degenerate
  forms $s(\tau^*\omega_2)_t+(1-s) (\omega_1)_t$, $s\in[0,1]$ provided
  $t\in (\lambda-\eps,\lambda+\eps)$ and $\eps<\eps_0$ for a
  sufficiently small $\eps_0$. Moser's lemma now provides a family of
  isotopies from $(\omega_1)_t$ to $(\tau^*\omega_2)_t$, which can be
  used to correct $\tau$ to the desired equivariant
  symplectomorphism. This shows the theorem in the case that the fixed
  point components are isolated points.

  The general case reduces to the case above as in \S\ref{sss:smooth}
  using normal directions. We have the following claim.

  Claim 2. The manifolds $M_i$ are isomorphic near the fixed
  point components.

  Suppose $\overline{F}_1\subset \red{M}{1,\lambda}$ is a fixed point
  component. The symplectomorphism $\phi$ of Definition \ref{d:sfpd1}
  maps $\overline{F}_1$ to a unique component $\overline{F}_2\subset
  \red{M}{2,\lambda}$. Moreover $\phi^*(N_{\overline{F}_2/
    \red{M}{2,\lambda}})=N_{\overline{F}_1/ \red{M}{1,\lambda}}$ and
  $\phi^*(P_{2,\lambda}) =P_{1,\lambda}$. Using Equations
  \eqref{eq:isobundle0}, \eqref{eq:isobundle1} and
  \eqref{eq:isobundle2} the normal bundles $N_{F_{i}/M}$ are also
  identified through $\phi$. Weinstein's symplectic neighborhood
  theorem shows that a neighborhood of $F_1$ is symplectomorphic to
  one of $F_2$. This proves Claim 2 and the Theorem.
\end{proof}

In \cite{Mc-So09} McDuff uses the previous argument to show that in the
case when the fixed point components are isolated points and the
reduced spaces are rigid (as in Definition \ref{d:rigid}) and
\emph{unique} in the sense that any two cohomologous symplectic forms
are diffeomorphic, then the germ is uniquely determined by the fixed
point data. We want to remark that it is not known whether rigidity
implies uniqueness.

Before we prove our next result, we need two more definitions. 

\begin{definition}\label{d:sfpd2}
  Let $(M_i,H_i,\omega_i), i=1,2$ denote two Hamiltonian
  $S^1$-manifolds.  We say that they have the \textbf{same small fixed
    point data} at a \textbf{non-extremal critical} value $\lambda$ if

\begin{enumerate} 
  \item $\lambda$ is a non-extremal critical value for both $H_i$. 
  \item There is a diffeomorphism 
    $\phi:(\red{M}{1,\lambda}, \overline{F}_{1,\lambda}) \to
    (\red{M}{2,\lambda}, \overline{F}_{2,\lambda})$, such that the
    restriction $\phi|_{F_1}$ is a symplectomorphism on each
    component.
  \item The index functions are related by $i_{2,\lambda}\circ \phi 
    =i_{1,\lambda}$.  
 \end{enumerate}
\end{definition}

Here we are only assuming the case when the fixed point components
have positive dimension in Definition \ref{d:mfpd}. The other cases
can be easily adapted. For the minimum, we assume the same as in
Definition \ref{d:sfpd1}. For the maximum we only require the
existence of a symplectomorphism $\phi:(F_{1,\max},
\omega_{F_{1,\max}})\to (F_{2,\max}, \omega_{F_{2,\max}})$ without
assuming anything about its normal bundle.

\begin{definition}
  Let $(M_i,H_i,\omega_i), i=1,2$ denote two Hamiltonian
  $S^1$-manifolds.  We say that the spaces $M_i$ are
  \textbf{identified} at a regular value $t$ if there is an
  $S^1$-bundle isomorphism $P_{1,t}\to P_{2,t}$ inducing a
  symplectomorphism of reduced spaces $\phi:(\red{M}{1,t},
  \overline{\omega}_{1,t})\to(\red{M}{2,t},
  \overline{\omega}_{2,t})$.
\end{definition}

What the next lemma says is that in the case that we have got two
manifolds with the same small fixed point data, then any
identification of the manifolds at a regular value $t$ right before
the critical level extends to an isomorphism over the critical level,
provided the family of reduced spaces is rigid. Note that this
identification plays the role of the symplectomorphism in Lemma
\ref{l:iso-germ}.

\begin{lemma}\label{l:iso-germ2}
  Let $(M_i,H_i,\omega_i), i=1,2$ denote two (possibly open)
  Hamiltonian $S^1$-manifolds with proper Hamiltonians $H_i$. Assume
  that $\lambda$ is a common critical value of $(M_i,H_i,\omega_i),
  i=1,2$ with only fixed point components of index 2. Suppose further
  that:
  
  \begin{enumerate}
  \item The manifolds have the same small fixed point data.
  \item $(M_i,H_i,\omega_i)$ are
    identified at a regular value $t_0<\lambda$.
  \item The pair
    ($\overline{M}_{1,t_0},\{\overline{\omega}_{1,t}\}_{t_0\leq
      t<\lambda}$) is rigid.
  \end{enumerate}
  Then, there is an $\eps_0>0$ such that $H_1^{-1}(\lambda-\eps,
  \lambda + \eps)$ is isomorphic to $H_2^{-1}(\lambda-\eps, \lambda +
  \eps)$, for all $\eps\leq\eps_0$.
\end{lemma}

\begin{proof}
  The proof is similar to that of Lemma \ref{l:iso-germ}. We only need
  two modifications.  First, the equivalent of Claim 2 follows since
  the manifolds are identified at the regular value $t_0$, thus they
  have isomorphic bundles $P_{t_0}$. If they have the same small fixed
  point data, the normal bundles $N_{\overline{F}_i/
    \red{M}{i,\lambda}}$ are also isomorphic. The same argument
  as in Claim 2 of Lemma \ref{l:iso-germ} shows that there are
  neighborhoods $U_i$ of $F_i$ and a symplectomorphism
  $\phi:(U_1,\omega_1)\to (U_2,\omega_2)$.

  We only need to show the equivalent of Claim 1, that is that the
  complements $H_i^{-1}(\lambda-\eps,\lambda+\eps)\backslash U_i$ are
  also symplectomorphic for $\eps$ small enough. Just as before, they
  are clearly diffeomorphic, since the action is free on this piece
  and by hypothesis there is an identification $\red{M}{1,t_0}\iso
  \red{M}{2,t_0}$ of the reduced spaces at regular level $t_0$.  Also,
  the family of reduced symplectic forms $\overline{\omega}_{1,t}$ is
  weakly equivalent to $\overline{\omega}_{2,t}$ for $t\in[t_0,t_1]$
  for any $t_1<\lambda$. By hypothesis ($\overline{M}_{1,t},
  \overline{\omega}_{1,t}$) spaces is rigid for $t\in [t_0,t_1]$. By
  taking $t_1$ sufficiently close to $\lambda$ and using the same
  arguments as in Lemma \ref{l:iso-germ} along with Theorem
  \ref{p:equiv} we conclude that $H_1^{-1}(\lambda-\eps,
  \lambda+\eps)\backslash U_1\iso H_2^{-1}(\lambda-\eps, \lambda +
  \eps ) \backslash U_2$.  Now, by taking $\eps$ small enough the
  symplectic structures on the two pieces $U_i$ and
  $H_i^{-1}(\lambda-\eps,\lambda+\eps)\backslash U_i $ are close
  enough. By Moser's Lemma we can assume they agree. This finishes the
  proof.
\end{proof}

\subsection{Non-simple levels.}\label{ss:non-simple}

We now describe how to deal with non-simple levels. Suppose
$(M,H,\omega)$ has a critical level at $\lambda$, with fixed point
components of index and coindex 2.  In Definition \ref{d:fpd-s}, the
index attached to $\lambda$ is not well defined, rather it is an
integer valued function $i_{\lambda}$ on the set of connected
components of the fixed point set. Therefore this definition needs to
be modified to accommodate this. Let $F_1$ denote the union of index 2
components, i.e. $F_1=i_{\lambda}^{-1}(2)$ and let $F_2$ denote the coindex
2 part.  The modification necessary to the make our previous results
work in this setting are as follows.

\begin{enumerate}
\item As before, we can decompose the normal bundles as positive and
  negative directions $N_{F_i/M}=N_i^- \oplus N^+_i, i=1,2$. A small
  neighborhood of the zero section in $N_{F_i/M}$ is isomorphic to
  a neighborhood $U_i=U_i^+\oplus U_i^-$ of $F_i$ in $M$. Let
  $\overline{X}_{1,t}$ denote the symplectic reduction $U_1^{-}\qu_{t}
  S^1$ if $t<\lambda$ or $U_1^{+}\qu_{t} S^1$ if $t>\lambda$, and let
  $\overline{X}_{2,t}$ denote the symplectic reduction $U_2^{+}\qu_{t}
  S^1$ if $t<\lambda$ or $U_2^{-}\qu_{t} S^1$ if $t>\lambda$.
\item The symplectic reduction $\red{M}{\lambda}$ at $\lambda$ also
  has a smooth structure. This is because the construction of the
  charts in \S \ref{sss:smooth} is done locally near $F_i$. We denote
  $\overline{F}_i$ the images of the fixed point set.
\item For $\eps>0$ the space $\red{M}{\lambda+\eps}$ is diffeomorphic
  to the blow-up of $\red{M}{\lambda}$ along $\overline{F}_1$
  only. $\red{M}{\lambda-\eps}$ is the blow up along
  $\overline{F}_2$. Let $\beta_{1}:\red{M}{\lambda+\eps}\to
  \red{M}{\lambda}$ and $\beta_{2}:\red{M}{\lambda-\eps}\to
  \red{M}{\lambda}$ denote the respective blow-down maps.  The tuple
  $(P_t, \red{M}{t},\overline{X}_{1,t}, \overline{X}_{2,t})$ depends
  smoothly on $t$, and its diffeotype is constant for $t$ close enough
  to $\lambda$. The main difference to the simple case is that none of
  reduced spaces for $\eps>0$ is diffeomorphic to the one at
  $\lambda$. Denote by $(P_{\pm}, \red{M}{\pm},\overline{X}_{1,\pm},
  \overline{X}_{2,\pm})$, the diffeotypes for $t=\lambda\pm \eps$
  respectively and $\eps>0$ small.  

\item There exist a smooth $S^1$ principal bundle $P$ over
  $\red{M}{\lambda}$ such that the principal bundles $P_{\pm}$ are
  related by
\begin{equation}\label{eq:bundles}
  \beta_1^* P=P_{+}\otimes L_1, \ \ \ \beta_2^*
  P=P_{-}\otimes L_2,
\end{equation}
for some circle bundles $L_i$ depending completely on the fixed point
data at $\lambda$, as in Equation \eqref{eq:euler}. More precisely,
Denote by $L_1'$ (resp. $L_2'$) the line bundle on $\red{M}{+}$
(resp. $\red{M}{-}$) whose Chern class is
dual to its codimension 2 submanifold  $\red{X}{1,+}$
(resp. $\red{X}{2,-}$). Let $L_1$ (resp. $L_2$) be the
circle bundle associated to $L_1'$ (resp. $L_2'$).

\item The family of $\eps$-blow up symplectic forms has a well defined
  germ, just as in \S\ref{sss:cobord}.

\item Definitions \ref{d:sfpd1} and \ref{d:sfpd2} are perfectly valid
  for the case at hand, since we require the morphism $\phi$ to
  preserves the index.

\item Claim 2 in the proofs of Lemma \ref{l:iso-germ} and Lemma
  \ref{l:iso-germ2} also holds in this case, since it is local in
  nature. Claim 1 also follows. Details are left to the reader.

\end{enumerate}

\section{Six dimensional case.}\label{s:examples}

In this section we will apply our previous analysis in dimension six.

\subsection{Some results in 4-dimensional symplectic topology.}

For our analysis, we will need the following results. The techniques
in our paper rely on them, and the possibility of extending our
results beyond six dimensions depends on the validity of
analogous theorems in higher dimensions.

\begin{thm}[\cite{Mc1}]\label{t:dusa}
  Let $X$ be $\CP{2}$, a blow-up of $\CP{2}$ or any rational
  surface. Then
  \begin{enumerate}
  \item Any deformation of two cohomologous symplectic forms on $X$
    may be homotoped through deformations with fixed endpoints to an
    isotopy. \footnote{This is actually true for more general
      4-manifolds, namely manifolds of non-simple SW-type.}
  \item Any two cohomologous symplectic forms on $X$ are
    symplectomorphic.
  \end{enumerate} 
\end{thm}

\begin{lemma}[Abreu, Gromov, McDuff, Lalonde, Pinsonnault,
  Evans]\label{l:idh}
  Let $(X,\om)$ be $\CP{2}$ with the form $\varpi_t(L)=t$ or its
  one-point blow-up $\CP{2} \# \overline{\CP{2}}$ such that on the
  exceptional divisor $\varpi_t(E)=t-\lambda_1$. Then the group
  $\mathrm{Symp} (X,\varpi_t)$ is connected for all $t$. Similarly if
  $(X,\om)$ denotes any of the blow ups $\CP{2} \# k\overline{\CP{2}}$
  for $k\le3$, $\varpi_t(E_i)=t-\lambda_i$ for the exceptional curves
  $E_i$, then the group $\mathrm{Symp}^{H} (X,\om)$ of
  symplectomorphisms that induce the identity on $H_{*}(X)$ is path
  connected.
\end{lemma}

The reader can consult the original articles \cite{Gr},\cite{AMc},
\cite{LaP}, \cite{P}, \cite{Ev-Sy09} the survey \cite{Mc3} for these
results. 

\subsection{Proofs of the main theorems.}\label{ss:A} 
\begin{figure}
\begin{center}
\setlength{\unitlength}{1900sp}%
\begingroup\makeatletter\ifx\SetFigFont\undefined%
\gdef\SetFigFont#1#2#3#4#5{%
  \reset@font\fontsize{#1}{#2pt}%
  \fontfamily{#3}\fontseries{#4}\fontshape{#5}%
  \selectfont}%
\fi\endgroup%
\begin{picture}(11709,4974)(3139,-5023)
\put(8551,-4561){\makebox(0,0)[lb]{\smash{{\SetFigFont{6}{7.2}
        {\rmdefault} 
        {\mddefault}{\updefault}{\color[rgb]{0,0,0} 
          H=0}%
}}}}
\thinlines
{\color[rgb]{0,0,1}
\put(9901,-3661){\line( 1, 2){450}}
{\color[rgb]{0,0,0}\put(10351,-2761){\line( 1, 0){900}}}
\put(11251,-2761){\line( 2,-3){450}}
\put(11701,-3436){\line(-4,-3){900}}
{\color[rgb]{0,0,0}\put(10801,-4111){\line(-2, 1){900}}}
}
{\color[rgb]{0,0,0}\put(12601,-511){\line( 1, 0){900}}
{\color[rgb]{0,0,1}\put(13501,-511){\line( 1,-1){450}}}
\put(13951,-961){\line( 0,-1){450}}
{\color[rgb]{0,0,1}\put(13951,-1411){\line(-2,-1){900}}}
\put(13051,-1861){\line(-2, 1){900}}
{\color[rgb]{0,0,1}\put(12151,-1411){\line( 1, 2){450}}}
}
{\color[rgb]{0,0,0}\put(9901,-511){\line( 0,-1){1350}}
\put(9901,-1861){\line( 1, 0){1350}}
\put(11251,-1861){\line(-1, 1){1350}}
}
{\color[rgb]{0,0,1}\put(12601,-2761){\line( 0,-1){1350}}
\put(12601,-4111){\line( 1, 0){1350}}
\put(13951,-4111){\line(-1, 1){1350}}
}

{\color[rgb]{1,0,0}\multiput(5626,-3886)(-368.18182,0.00000){6}{\line(-1, 0){184.091}}
\multiput(3601,-3886)(959.89387,-319.96462){2}{\line( 3,-1){390.106}}
\put(4951,-4336){\line( 3, 2){675}}
}

{\color[rgb]{.7,.5,.9}
\multiput(3151,-2311)(19.18605,0.00000){259}{\makebox(3.1750,22.2250){\SetFigFont{5}{6}{\rmdefault}{\mddefault}{\updefault}.}}
\multiput(3151,-3661)(19.18605,0.00000){259}{\makebox(3.1750,22.2250){\SetFigFont{5}{6}{\rmdefault}{\mddefault}{\updefault}.}}
\multiput(4951,-3211)(19.20732,0.00000){165}{\makebox(3.1750,22.2250){\SetFigFont{5}{6}{\rmdefault}{\mddefault}{\updefault}.}}
\multiput(4951,-4561)(19.20732,0.00000){165}{\makebox(3.1750,22.2250){\SetFigFont{5}{6}{\rmdefault}{\mddefault}{\updefault}.}}
\multiput(6751,-2761)(19.28571,0.00000){71}{\makebox(3.1750,22.2250){\SetFigFont{5}{6}{\rmdefault}{\mddefault}{\updefault}.}}
\multiput(6751,-1411)(19.28571,0.00000){71}{\makebox(3.1750,22.2250){\SetFigFont{5}{6}{\rmdefault}{\mddefault}{\updefault}.}}
\multiput(4951,-1861)(19.20732,0.00000){165}{\makebox(3.1750,22.2250){\SetFigFont{5}{6}{\rmdefault}{\mddefault}{\updefault}.}}
\multiput(4951,-511)(19.20732,0.00000){165}{\makebox(3.1750,22.2250){\SetFigFont{5}{6}{\rmdefault}{\mddefault}{\updefault}.}}

}

\thicklines
{\color[rgb]{0,0,0}\put(8000,-61){\line( 0,-1){4950}}
}
{\color[rgb]{0,0,0}\multiput(4951,-1861)(716.10221,-358.05110){3}{\line( 2,-1){367.796}}
}
{\color[rgb]{0,0,0}\multiput(4951,-511)(0.00000,-385.71429){4}{\line( 0,-1){192.857}}
\multiput(4951,-511)(-503.07735,-503.07735){4}{\line(-1,-1){290.768}}
}
{\color[rgb]{0,0,0}\multiput(3151,-3661)(500,500){2}{\line(1,1){250}}
}


{\color[rgb]{0,0,0}\put(3151,-2311){\line( 1, 1){1800}}
}
{\color[rgb]{0,0,0}\put(6751,-1411){\line(-2, 1){1800}}
}
{\color[rgb]{0,0,0}\put(4951,-3211){\line( 1, 1){1800}}
}
{\color[rgb]{0,0,0}\put(4951,-3211){\line(-2, 1){1800}}
}
{\color[rgb]{0,0,0}\put(3151,-3661){\line( 0, 1){1350}}
}
{\color[rgb]{0,0,0}\put(6751,-2761){\line( 0, 1){1350}}
}
{\color[rgb]{0,0,0}\put(4951,-4561){\line( 1, 1){1800}}
}
{\color[rgb]{0,0,0}\put(4951,-4561){\line( 0, 1){1350}}
}
{\color[rgb]{0,0,0}\put(4951,-4561){\line(-2, 1){1800}}
}
\thinlines

\put(13051,-4786){\makebox(0,0)[lb]{\smash{{\SetFigFont{6}{7.2}{\rmdefault}{\mddefault}{\updefault}{\color[rgb]{0,0,0}$d)$}%
}}}}
\put(10351,-4786){\makebox(0,0)[lb]{\smash{{\SetFigFont{6}{7.2}{\rmdefault}{\mddefault}{\updefault}{\color[rgb]{0,0,0}$c)$}%
}}}}
\put(13051,-2311){\makebox(0,0)[lb]{\smash{{\SetFigFont{6}{7.2}{\rmdefault}{\mddefault}{\updefault}{\color[rgb]{0,0,0}$b)$}%
}}}}
\put(10351,-2311){\makebox(0,0)[lb]{\smash{{\SetFigFont{6}{7.2}{\rmdefault}{\mddefault}{\updefault}{\color[rgb]{0,0,0}$a)$}%
}}}}
\put(9500,-4361){\makebox(0,0)[lb]{\smash{{\SetFigFont{6}{7.2}{\rmdefault}{\mddefault}{\updefault}{\color[rgb]{0,0,0}$L-E_{1}-E_{2} = 0$}%
}}}}
\put(11300,-1861){\makebox(0,0)[lb]{\smash{{\SetFigFont{6}{7.2}{\rmdefault}{\mddefault}{\updefault}{\color[rgb]{0,0,0}$L-E_{2}-E_{3}$}%
}}}}
\put(13440,-1186){\makebox(0,0)[lb]{\smash{{\SetFigFont{6}{7.2}{\rmdefault}{\mddefault}{\updefault}{\color[rgb]{0,0,0}$L-E_{1}-E_{2}$}%
}}}}
\put(12200,-400){\makebox(0,0)[lb]{\smash{{\SetFigFont{6}{7.2}{\rmdefault}{\mddefault}{\updefault}{\color[rgb]{0,0,0}$L-E_{1}-E_{3}$}%
}}}}
\put(12000,-736){\makebox(0,0)[lb]{\smash{{\SetFigFont{6}{7.2}{\rmdefault}{\mddefault}{\updefault}{\color[rgb]{0,0,1}$E_{3}$}%
}}}}
\put(13501,-1861){\makebox(0,0)[lb]{\smash{{\SetFigFont{6}{7.2}{\rmdefault}{\mddefault}{\updefault}{\color[rgb]{0,0,1}$E_{2}$}%
}}}}
\put(13951,-736){\makebox(0,0)[lb]{\smash{{\SetFigFont{6}{7.2}{\rmdefault}{\mddefault}{\updefault}{\color[rgb]{0,0,1}$E_{1}$}%
}}}}
\put(10351,-2086){\makebox(0,0)[lb]{\smash{{\SetFigFont{6}{7.2}{\rmdefault}{\mddefault}{\updefault}{\color[rgb]{0,0,0}$L$}%
}}}}
\put(10801,-1186){\makebox(0,0)[lb]{\smash{{\SetFigFont{6}{7.2}{\rmdefault}{\mddefault}{\updefault}{\color[rgb]{0,0,0}$L$}%
}}}}
\put(9676,-1186){\makebox(0,0)[lb]{\smash{{\SetFigFont{6}{7.2}{\rmdefault}{\mddefault}{\updefault}{\color[rgb]{0,0,0}$L$}%
}}}}
\put(8200,-511){\makebox(0,0)[lb]{\smash{{\SetFigFont{6}{7.2}{\rmdefault}{\mddefault}{\updefault}{\color[rgb]{0,0,0}$\la_{1} +\la_{2}+ \la_{3}$}%
}}}}
\put(8200,-1411){\makebox(0,0)[lb]{\smash{{\SetFigFont{6}{7.2}{\rmdefault}{\mddefault}{\updefault}{\color[rgb]{0,0,0}$\la_{2}+\la_{3}$}%
}}}}
\put(8200,-1861){\makebox(0,0)[lb]{\smash{{\SetFigFont{6}{7.2}{\rmdefault}{\mddefault}{\updefault}{\color[rgb]{0,0,0}$\la_{1}+\la_{3}$}%
}}}}
\put(8200,-2311){\makebox(0,0)[lb]{\smash{{\SetFigFont{6}{7.2}{\rmdefault}{\mddefault}{\updefault}{\color[rgb]{0,0,0}$\la_{1}+\la_{2}$}%
}}}}
\put(8200,-2761){\makebox(0,0)[lb]{\smash{{\SetFigFont{6}{7.2}{\rmdefault}{\mddefault}{\updefault}{\color[rgb]{0,0,0}$\la_{3}$}%
}}}}
\put(8200,-3211){\makebox(0,0)[lb]{\smash{{\SetFigFont{6}{7.2}{\rmdefault}{\mddefault}{\updefault}{\color[rgb]{0,0,0}$\la_{2}$}%
}}}}
\put(8200,-3661){\makebox(0,0)[lb]{\smash{{\SetFigFont{6}{7.2}{\rmdefault}{\mddefault}{\updefault}{\color[rgb]{0,0,0}$\la_{1}$}%
}}}}

\end{picture}%
\end{center}
\caption{ The fixed points,  their critical values and some reduced
  spaces of a manifold with isolated fixed points (Example A).
  Figures $a), b), c)$ and $d)$ represent the reduced spaces for $t$
  in the intervals
  $(0,\la_{1}),(\la_{3},\la_{1}+\la_{2}),(\la_{1}+\la_{2},\la_{1}+\la_{3}),
  (\la_{2}+\la_{3},\la_{1}+\la_{2}+\la_{3})$ respectively. Note that
  as $t\to \la_{1}+\la_{2}$ the exceptional sphere $L-E_{1}-E_{2}$
  blows down.  $d)$ is the manifold obtained after blowing down the
  exceptional spheres $L-E_{i}-E_{j}$. $a)$ and $d)$ are diffeomorphic
  (via  Cremona transformation) but the Euler class of the principal
  bundles associated to these reduced spaces differ by a sign.} 
\label{f:2}
\end{figure}
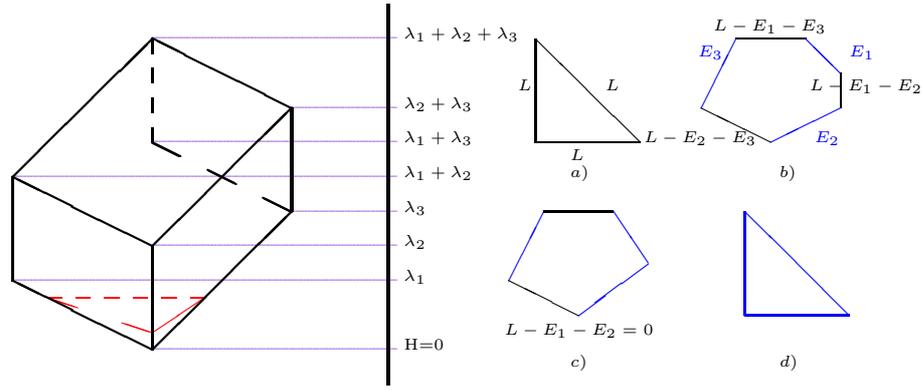

Assume that $(M,\om)$ is a symplectic 6-manifold with a semi-free
$S^{1}$-action with isolated fixed points. Since the fixed points are
isolated it is not necessary to assume that the action is Hamiltonian,
it would follow from \cite{TW}.  The fixed points are given as
follows.  The minimum, 
three critical points $p_{1}, p_{2}, p_{3}$ of index two, three
$p_{12}, p_{23}, p_{13}$ of index four and a maximum. We will denote by
$\la_{i}$ the critical values $H(p_{i})$.  Without lost of generality
we may assume the minimum value of $H$ is zero, and that
$\la_{1}\le\la_{2}\le\la_{3}$.  It is not hard to see \cite{G1} that
the fixed points of index 4 are in the level sets
$H^{-1}(\la_{i}+\la_{j}), j\ne i$ and the maximum is the unique point
in $H^{-1}(\la_{1}+\la_{2}+\la_{3})$ (See Figure \ref{f:2}).  Before
going any further recall that we want to prove that $M$ is isomorphic
to $Y^{3}=S^{2}\times S^{2}\times S^{2}$ with the product symplectic
form $\varpi=\la_{1}\sg\times\la_{2}\sg\times \la_{3}\sg$.  Here $\sg$
is the canonical area form on $S^{2}$.  We are assuming that the
circle acts on $Y^{3}$ by
\[
e^{2\pi i t}(x,y,z) \mapsto (e^{2\pi it}x, e^{2\pi it}y, e^{2\pi it}z).
\]
and denote by $K:Y\to \bb{R}$ its Hamiltonian. For simplicity of the
notation we write $Y^t:=K^{-1}[0,t)$ and similarly $M^t:=H^{-1}[0,t)$.

\begin{proof}[Proof of Theorem \ref{t:intro}]
  We start by noticing that there are two main cases to analyse (but
  not only), when $\la_{1}+\la_{2}\leq \la_{3}$ and when
  $\la_{1}+\la_{2}>\la_{3}$.  The difference of this two cases is the
  order in which we reach fixed points.  For simplicity, one can treat
  the first case since the second one is analogous. We start by
  assuming that none of the $\lambda_i$ agree, since this case is
  slightly easier than the case where some of them are equal.

  Our aim is to give a construction of the isomorphism by showing that
  the manifolds $M^t$ and $Y^t$ are isomorphic as $t$ increases. We
  will start from the minimum level.  We have already explain that the
  minimal germs are determined by the fixed point data at the minimum
  (see comment after Definition \ref{d:1}).  We now explain this for
  this particular example, since we want to be as explicit as
  possible.  From the equivariant version of the Darboux theorem at
  the minimum, one gets a neighborhood of the minimum isomorphic to
  $\bb{C}^{3}$ with the diagonal circle action.  The Hamiltonian
  function in these coordinates is given by
  \[
  (z_{1}, z_{2}, z_{3})\mapsto |z_{1}|^{2} +|z_{2}|^{2} +|z_{3}|^{2}.
  \]
  If $0<t<\eg$ for $\eg>0$ small enough, we have that the level set
  $H=t$ is the 5-sphere,
  \[
  S^{5} \ : \  |z_{1}|^{2} +|z_{2}|^{2} +|z_{3}|^{2}=t.
  \]
  Therefore, the symplectic reduction $\red{M}{t}=S^{5}/S^{1}$ is the
  projective space $(\CP{2},\overline{\omega}_{t})$ with the
  symplectic form $\overline{\omega}_{t}$ that takes the value $t$ on
  the line $L$.  The bundle $P=H^{-1}(t)\to \red{M}{t}$ is now just
  the Hopf fibration $S^{5}\to \CP{2}$, whose Euler class $e(P)\in
  H^{2}(\CP{2})$ is the negative generator of $H^{2}(\CP{2})$.  The
  diffeomorphism type of the reduced space $Y_{t}$ for all
  $0<t<\lambda_{1}$ is $\CP{2}$.  Theorem \ref{t:dusa} and Lemma
  \ref{l:idh} show that the pair
  $(\CP{2},\{\overline{\omega}_t\}_{t\in I})$ is rigid for all
  subinterval $I\subset (0,\lambda_1)$, and thus by Theorem
  \ref{p:equiv} the symplectomorphism type of the slice corresponding
  to the interval $I$ is determined. The diffeotype of the reduced
  spaces and the respective bundles is constant for all $t\leq
  \lambda$, in particular for $t=\lambda$. Thus (Definition
  \ref{d:sfpd2}) both manifolds have the same fixed point data at
  $\lambda_1$ and they are identified at a regular level immediately
  before $\lambda_1$.  Lemma \ref{l:iso-germ2} shows we can extend the
  isomorphism pass the critical level. One then obtains an isomorphism
  $M^{\lambda_1+\eps}\iso Y^{\lambda_1+\eps}$ by gluing along any
  regular level $t<\lambda_1$. It is important to note that \emph{any}
  gluing map $\Phi:H^{-1}(t)\to K^{-1}(t)$ would suffice, since we
  just want to show the existence of the isomorphism.
 
  For $t-\lambda_1>0$ small enough, Theorem \ref{t:bordism} asserts
  that the reduction at $t$ is the blowup $\CP{2}\#\overline{\CP{2}}$
  with the symplectic form $\overline{\omega}_t$ such that the line
  class $L$ has symplectic area $t$ and the exceptional class $E_1$
  has symplectic area $t-\lambda_1$. The principal bundle $P_t$
  corresponding to the level $t$ is explicitly given in terms of the
  previous fixed point data Example \ref{ex:1}. The germ near
  $\lambda_1$ is determined by the fixed point data at $t=\lambda_1$,
  so we have $M^t\iso Y^t$ for $t-\lambda_1>0$. Since $(\CP{2}\#
  \overline{\CP{2}}, \{\overline{\omega}_t\}_{t\in (\lambda_1,\lambda_2)})$ is
  rigid, the same argument above, shows that $M^{\lambda_2}\iso
  Y^{\lambda_2}$.

  The symplectic area of the exceptional class $E_1$ in $\CP{2}\#
  \overline{\CP{2}}$ depends linearly in $t$, and thus cannot
  blow-down as $t$ reaches $\lambda_{2}$.  For $t-\lambda_2>0$ small
  enough the reduced space blows up and thus it is
  $\CP{2}\#2\overline{\CP{2}}$ equipped with the family of blow up
  symplectic forms $\{\overline{\omega}_t\}$ taking values $t$,
  $t-\lambda_1$, $t-\lambda_2$ on $L,E_1,E_2$ respectively. The pair
  $(\CP{2}\#2\overline{\CP{2}},\{\overline{\omega}_t\}_{t\in(\lambda_1,\lambda_2)})$
  is again by Theorem \ref{t:dusa} and Lemma \ref{l:idh} rigid, and
  then we have got an isomorphism $Y^t\iso M^t$ for
  $t<\lambda_1+\lambda_2$.  To see what happens after we pass the
  critical level $\lambda_1+\lambda_2$, we note that there are only
  three exceptional classes: $E_1, E_2$ and $E_{12}:=L-E_1-E_2$ with
  respective areas given by $t-\lambda_1$, $t-\lambda_2$ and
  $t-\lambda_1-\lambda_2$. The areas of the exceptional curves
  $E_1,E_2$ are growing, thus do not blow down as $t$ reaches
  $\lambda_1+\lambda_2$, then only a curve in class $E_{12}$ blows
  down when crossing $\lambda_1+\lambda_2$.  One important issue is to
  note that there is only one way of blowing-down. This is, regardless
  of the actual curve that blows down, the symplectomorphic type of
  the blow-down manifold only depends on the homology class of the
  curves \cite[proof of Theorem 1.1]{Mc:Ra}. In our case the class
  $E_{12}$ is the one that blows down, and thus we have got an
  isomorphism $Y^t\iso M^t$ for $t=\lambda_1+\lambda_2+\eps$ and
  $\eps>0$ small enough.  Following this process, one notices that the
  only possible reduced spaces that can appear are $\CP{2} \#
  k\overline{\CP{2}}$ for $k\leq 3$, which have associated rigid
  families of symplectic structures, by Theorem \ref{t:dusa} and Lemma
  \ref{l:idh}. The classes that blow down are determined by the
  numbers $\lambda_i$ and thus the same arguments we have used before
  give the desired isomorphism $M^t\iso Y^t$.  The process continues
  until $t$ reaches the maximum, obtaining the desired isomorphism, in
  the case when all the critical values $\lambda_i$ are different.

  The case when two or three of the values $\lambda_i$ agree, one has
  to pay attention to the validity of condition \ref{d:pathconnected}
  in Definition \ref{d:rigid}; there might be several components in
  the group of symplectomorphisms induced by symplectomorphisms that
  permute exceptional classes of equal size in the blow
  up. Nevertheless, the reduced spaces are still rigid, since the
  symplectomorphisms in the base that we are considering are those
  which also are in the identity component of the diffeomorphisms
  which is connected by \cite{P,Ev-Sy09}.

\end{proof}

The more general results stated in the introduction now follow
similarly.

\begin{proof}[Proof of Theorem \ref{t:intro1}]
  Recall that $(M,H,\om)\in \HS$ and $\Cc(M)$ is the set of critical
  values.  We are assuming that each non-extremal critical value
  $\la\in \Cc(M)=\{\la_0,\dots, \la_{s}\}$ contains only fixed points
  of (co)index 2 and that all the families of reduced spaces
  $(\red{M}{\la},
  \overline{\omega}_t)_{t\in(\lambda_i,\lambda_{i+1})}$ for
  consecutive critical values are rigid. Therefore, by Lemma
  \ref{l:iso-germ} and Theorem \ref{p:equiv} the regular slices and
  germs are completely determined by the fixed point data. Given any
  other manifold with the same fixed point data. We use the same
  procedure as in the proof of Theorem \ref{t:intro}, to build the
  isomorphism starting from the minimum and gluing slices and germs. 
\end{proof}

\begin{proof}[Proof of Theorem \ref{t:baby2}]
  Recall that we are assuming that for each $\la\in\Cc(M)$ the fixed
  point submanifolds
  $(\overline{F}_{\la_i},\overline{\omega}_{F_{\la_i}})$ are isolated
  points or surfaces of index 2. 

  The proof is again by exhaustion on the sets
  $H^{-1}(-\infty,t)$. Following the same notations as before, suppose
  $\Cc(M)=\{\lambda_0=0,\lambda_1, \dots , \lambda_s\}$.  Starting at
  the minimum, the normal form given by the Darboux-Weinstein Theorem
  applied to the fixed point set at $\lambda_0$, shows that the
  minimal germ is determined by the index and the submanifold
  $(\overline{F}_{\la_0}, \overline{\omega}_{F_{\la_0}})$ (with its
  normal bundle). By taking the reduction at any level $t$ close
  enough to $\la_0$ one obtains the diffeomorphism type of the reduced
  bundle $P_{t}\to \red{M}{t}$.  The rigidity of the reduced spaces
  and Theorem \ref{p:equiv} ensures that the isomorphism type of the
  regular slice for the interval $(\lambda_0,\lambda_1)$ is determined
  by the small fixed point data at $\lambda_0,\lambda_1$.

  Let $t=\lambda_1-\eps$, for $\eps>0$ small enough, thus any other
  manifold with the same small fixed point data is isomorphic to $M$
  up to level $t$ and thus they are identified at level $t$
  (cf. Definition \ref{d:sfpd2}), satisfying the hypotheses of Lemma
  \ref{l:iso-germ2}. This shows that the critical germ at $\lambda_1$
  is determined, and therefore any other space with the same fixed
  point data is identified with $M$ up to level $t=\lambda_1+\eps$,
  for some $\eps>0$ small enough.  Examples \ref{ex:1}, \ref{ex:2}
  show that the principal bundle $P_{\la_2}\to \red{M}{\la_2}$ is also
  determined by the known data at $\la_1$. We continue this process,
  until $t$ reaches $\lambda_s$.
\end{proof}

Finally, using the analysis of Section \ref{ss:non-simple} one can show
that this theorems hold even in the non-simple case. The details are
left to the reader.

\bibliographystyle{plainnat}
\bibliography{diss}
\end{document}